\documentclass[aos,preprint]{imsart}

\usepackage{hyperref}
\usepackage{graphicx}
\usepackage{tikz, pgfplots}
\usepackage{subcaption}
\usepackage{bm}
\definecolor{dodgerblue}{rgb}{0.12, 0.56, 1.0}
\definecolor{firebrick}{rgb}{0.7, 0.13, 0.13}
\definecolor{goldenrod}{rgb}{0.85, 0.65, 0.13}

\newcommand{\p}[1]{\left(#1\right)}

\newcommand{\argmax}{\operatorname{argmax}}

\newcommand{\sign}{\operatorname{sign}}

\newcommand{\hf}{\hat{f}}

\newcommand{\eqif}{\text{ if }}

\usepackage{amsthm,amsmath,amssymb}
\usepackage[authoryear]{natbib}
\usepackage{caption, subcaption}
\usepackage{graphicx} 
\usepackage{algorithm}
\usepackage{comment}
\usepackage{algorithm}
\usepackage[noend]{algpseudocode}

\newcommand{\E}{\mathbb{E}}
\renewcommand{\P}{\mathbb{P}}
\newcommand{\R}{\mathbb{R}}
\newcommand{\A}{\mathcal{A}}

\newcommand{\supp}{\textnormal{supp}}
\newcommand{\bitem}{\begin{itemize}}
\newcommand{\eitem}{\end{itemize}}
\newcommand{\benum}{\begin{enumerate}}
\newcommand{\eenum}{\end{enumerate}}
\newcommand{\beq}{\begin{equation}}
\newcommand{\eeq}{\end{equation}}
\newcommand{\beqs}{\begin{equation*}}
\newcommand{\eeqs}{\end{equation*}}

\def\balign#1\ealign{\begin{align}#1\end{align}}
\def\baligns#1\ealigns{\begin{align*}#1\end{align*}}
\def\reals{\mathbb{R}} 

\newcommand{\ind}[1]{\mathbb{I}_{\{#1\}}}

\newcommand{\iid}{\stackrel{\mathrm{iid}}{\sim}}

\newtheorem{theorem}{Theorem}[section]
\newtheorem{definition}{Definition}
\newtheorem{lemma}{Lemma}[section]
\newtheorem{claim}{Claim}[section]

\newtheorem{rmk}{Remark}
\theoremstyle{thm}

\arxiv{arXiv:0000.0000}

\begin{document} 

\begin{frontmatter}
\title{Transfer Learning for Nonparametric Classification: Minimax Rate and Adaptive Classifier\thanksref{T1}}
\runtitle{Transfer Learning}
\thankstext{T1}{The research was supported in part by NSF Grant DMS-1712735 and NIH grants R01-GM129781 and R01-GM123056.}
\begin{aug}
\author{\fnms{T. Tony} \snm{Cai}\ead[label=e1]{tcai@wharton.upenn.edu}}
\and
\author{\fnms{Hongji} \snm{Wei}\ead[label=e2]{hongjiw@wharton.upenn.edu}
\ead[label=u1,url]{http://www-stat.wharton.upenn.edu/$\sim$tcai/}}
\runauthor{T. T. Cai and H. Wei}
\affiliation{University of Pennsylvania}
\address{DEPARTMENT OF STATISTICS\\
THE WHARTON SCHOOL\\
UNIVERSITY OF PENNSYLVANIA\\
PHILADELPHIA, PENNSYLVANIA 19104\\
USA\\
\printead{e1}\\
\phantom{E-mail:\ }\printead*{e2}\\
\printead{u1}\\
}

\end{aug}

\begin{abstract}
Human learners have the natural ability to use knowledge gained in one setting for learning in a different but related setting. This ability to transfer knowledge from one task to another is essential for effective learning. In this paper, we study transfer learning  in the context of nonparametric classification based on observations from different distributions under the posterior drift model, which is a general framework and arises in many practical problems. 

We first establish the minimax rate of convergence and construct a rate-optimal two-sample weighted $K$-NN classifier. The results characterize precisely the contribution of the observations from the source distribution to the classification task under the target distribution. A data-driven adaptive classifier is then proposed and is shown to simultaneously attain within a logarithmic factor of the optimal rate over a large collection of parameter spaces. Simulation studies and real data applications are carried out where the numerical results further illustrate the theoretical analysis.  Extensions to the case of multiple source distributions are also considered. 
\end{abstract}



\end{frontmatter}
	
\section{Introduction}
	
	A key feature of intelligence is the ability to learn from experience.  Human learners appear to have the talent to transfer their knowledge gained from one task to another similar but different task. However, in statistical learning, most procedures are designed to solve one single task, or to learn one single distribution based on observations from the same setting. In a wide range of real-world applications, it is important to gain improvement of learning in a new task through the transfer of knowledge from a related task that has already been learned. Transfer learning aims to tackle such a problem. It has attracted increasing attention in machine learning and has been used in many applications. Recent examples include computer vision \citep{tzeng2017adversarial, gong2012geodesic}, speech recognition \citep{huang2013cross}, genre classification \citep{choi2017transfer} and also many newly designed algorithms such as \cite{yao2010boosting, lee2007learning}. More details about transfer learning can be found in the survey papers \citep{pan2010survey, weiss2016survey}.
	 
	Besides significant successes in applications, much recent focus has also been on the theoretical properties of transfer learning. In many practical situations, there are labeled data available from a distribution $P$, called the source distribution, while a relatively small quantity of labeled or unlabeled data is drawn from a distribution $Q$, called the target distribution. They are different but to some extent related distributions. The goal is to make statistical inference under $Q$. A natural questions is: How much information can be transferred from the source distribution $P$ to the target distribution $Q$, provided a certain level of similarity between the two distributions?
 	
	This is quite a general and challenging question. The problem is also known as domain adaptation in the binary classification setting. In domain adaptation, data pairs $(X,Y)$ are drawn from $P$ and $Q$ defined on $\reals^d \times \{0,1\}$. Data from the source distribution $P$ can be informative about the target distribution $Q$ if the two distributions are similar.  Several type of assumptions have been proposed and studied previously in the literature, such as divergence bounds, covariate shift, and posterior drift. 
	The first line of work in the literature measures the similarity by the divergence between $P$ and $Q$. Generalization bounds are derived on unlabeled testing data from the target distribution $Q$ after training by the data from the source distribution $P$ \citep{ben2007analysis, blitzer2008learning, mansour2009domain}. These bounds are general and can be applied to any two distributions, but for more structured source and target distributions those bounds are not suitable.	
	Another line of work in the literature imposes some structural assumptions on $P$ and $Q$ such as covariate shift and posterior drift. Covariate shift assumes that the conditional distributions of $Y$ given $X$ are the same under $P$ and $Q$, i.e. $P_{Y|X} = Q_{Y|X}$, but the marginal distributions $P_X$ and $Q_X$ can be different. Such kind of setting typically arises when the same study/survey is carried out in different populations. For example, when constructing a classifier for a certain disease, source data may be generated from clinical studies, but  the goal is to classify people drawn from the general public. The task can become challenging due to the difference between the two populations. Transfer learning under the  covariate shift framework has been studied in previous works such as \cite{shimodaira2000improving, sugiyama2008direct, kpotufe2018marginal}.
	 
	In the present paper, we study transfer learning under the posterior drift model, where it is assumed that  $P_X = Q_X$ but $P_{Y|X}$ and $Q_{Y|X}$ can be different. To be more specific, suppose there are two data generating distributions $P$ and $Q$ on $[0,1]^d \times \{0,1\}$. We observe $n_P$ independent and identically distributed (i.i.d.) samples $(X_1^P, Y_1^P),  ..., (X_{n_P}^P, Y_{n_P}^P)$ drawn from a source distribution $P$, and $n_Q$ i.i.d. samples $(X_1^Q, Y_1^Q), ..., (X_{n_Q}^Q, Y_{n_Q}^Q)$ drawn from a target distribution $Q$. The data points from the distributions $P$ and $Q$ are also mutually independent. In each pair of data $(X,Y)$, the $d$-dimensional vector $X$ is regarded as covariates (features) of a certain object, while $Y$ is a (noisy) binary label indicating which of the two classes this object belongs to. The goal is to make classification under the target distribution $Q$: Given the observed data, construct a classifier $\hat f: [0,1]^d \to \{0,1\}$ which minimizes the classification risk under the target distribution $Q$:
	\[   R(\hf) \triangleq \P_{(X,Y)\sim Q}(Y \neq \hf(X) ),    \]
	here $\P_{(X,Y)\sim Q}(\cdot)$ means the probability when $(X,Y)$ are drawn from distribution $Q$.
	
	In such binary classification problems,  the regression functions are defined as
	\[ \eta_P(x) \triangleq P(Y=1|X=x) \quad{\rm and}\quad  \eta_Q(x) \triangleq Q(Y=1|X=x),  \]
	which can be used to represent the conditional distributions $P_{Y|X}$ and $Q_{Y|X}$. In classification, $Y$ can be regarded as an unknown parameter predicted by $X$, so from this perspective we refer to $P_X$ and $Q_X$ as the class ``prior" probabilities and $\eta_P(x)$ and $ \eta_Q(x)$ as the class ``posterior" probabilities associated with $P$ and $Q$ respectively \citep{scott2018generalized}. We say a ``posterior drift" happens when $P_X = Q_X$ but $\eta_P(x) \neq \eta_Q(x)$.

	Posterior drift is a general framework and broadly arises in many practical problems, where one collects data from different populations. Here are three examples.
	
	\begin{itemize}
		\item {\bf Crowdsourcing.} Crowdsourcing is a distributed model for large-scale problem-solving and experimentation such as image classification, video annotation, and translation  \citep{yuen2011survey, karger2011budget, zhang2014spectral}. The tasks are broadcasted to multiple independent workers online in order to collect and aggregate their solutions. In crowdsourcing, many noisy answers/labels are available from a large amount of public workers, while sometimes, more accurate answers/labels may be collected from experienced workers or experts. These expert answers/labels are of higher quality but are relatively few due to the time or budget constraints. One can view this difference in labeling accuracy as a posterior drift. It is desirable to construct a statistical procedure that incorporates both data sets.
		
		\item {\bf Concept drift.} Concept drift is a common phenomenon when the underlying distribution of the data changes over time in a streaming environment \citep{tsymbal2004problem, gama2014survey}. One kind of concept drift is called real concept drift where the posterior class probabilities $P(Y|X)$ changes over time. In this situation, posterior drift exists if data are collected at different time. For example, the incidence rate of a certain disease in certain groups may change over time due to the development of treatments and preventive measures.

		\item {\bf Data corruption.} Data corruption is ubiquitous in statistical applications, where unexpected error on data occurs during storage, transmission or processing \citep{menon2015learning, rooyen2018a}. In many settings, one receives data of variable quality -- perhaps some small amount of entirely clean data, another amount of slightly corrupted data, yet more that is significantly corrupted, and so on \citep{crammer2006learning}. Data of variable qualities can be viewed as posterior drift between those data generating distributions, thus better strategies are needed to tackle with the posterior drift between those data of variable qualities.
		
		\end{itemize}
				
	Under the posterior drift model, the difference between $P$ and $Q$ lies in the regression functions $\eta_P(x)$ and $\eta_Q(x)$. So the relationship between $\eta_P(x)$ and $\eta_Q(x)$, which can be captured by the link function $\phi$ defined below, is important in  characterizing the difficulty of the transfer learning problem. In this work, we propose a new concept called the {\it relative signal exponent} $\gamma$ to describe the relationship between $\eta_P(x)$ and $\eta_Q(x)$. Our results show that the relative signal exponent $\gamma$ plays an important role in the minimax rate of convergence for the excess risk under the  posterior drift model.
	
	For conceptual simplicity, we assume $\eta_P(x) = \phi(\eta_Q(x))$ for some strictly increasing link function $\phi(\cdot)$ with $\phi(\frac{1}{2}) = \frac{1}{2}$. Note that this is only a simplified version of our formal model which will be given in Section 2. It is natural to assume $\phi$ is strictly increasing in the settings where those $X$ that are more likely to be labeled $Y=1$ under $Q$ are also more likely to be labeled $Y=1$ under $P$. The assumption $\phi(\frac{1}{2}) = \frac{1}{2}$ means that those $X$ that are non-informative under $Q$ are the same under $P$. For a  given relative signal exponent $\gamma > 0$ and a constant $C_\gamma>0$, we denote by $\Gamma(\gamma, C_\gamma)$ the collection of all distribution pairs $(P,Q)$ satisfying
	\beq 
	 (\phi(x)-\frac{1}{2})(x-\frac{1}{2}) \geq 0 \quad{\rm and}\quad  |\phi(x)-\frac{1}{2}| \geq C_\gamma |x-\frac{1}{2}|^\gamma.  \label{intro1}  
	\eeq
	
	The relative signal exponent is a key parameter in capturing the usefulness of the data from the source distribution $P$ for the task of classification under $Q$. The smaller the relative signal exponent, the more information transferable from the source distribution $P$ to the target distribution $Q$, and vice versa. 
	
	In this work we consider transfer learning under the posterior drift model in a nonparametric classification setting. When $Q$ satisfies the margin assumption with the parameter $\alpha$, defined in Section \ref{Model.sec}, and $\eta_Q(x)$ belongs to the $(\beta, C_\beta)$-H\"{o}lder function class, it is shown that, under the regularity conditions,  the minimax optimal rate of convergence is given by
	\beq 
	\inf_{\hat f} \max_{(P,Q)\in \Pi} \E_Z \mathcal{E}_Q(\hat f) \asymp (n_P^{\frac{2\beta + d}{2\gamma\beta+d}} + n_Q)^{-\frac{\beta(1+\alpha)}{2\beta+d}}, \label{intro_rate} 
	\eeq
	where $n_P$ and $n_Q$ are number of data drawn from $P$ and $Q$ respectively, $d$ is the number of features, and $\Pi$ is the posterior drift regime where the distribution pair $(P,Q)$ belongs to the class $\Gamma(\gamma, C_\gamma)$ with the relative signal exponent $\gamma$ and satisfies some additional regularity conditions. Here $\mathcal{E}_Q(\hf)$ is the excess risk (also called regret) on $Q$ which is defined based on the misclassification error:
	\beq  \mathcal{E}_Q(\hf) = R_Q(\hf) - R_Q(f_Q^\ast) \label{def_exrisk} \eeq
	where 
	\beq f^\ast_Q(x) = \begin{cases} 0 \text{ if } \eta_Q(x) \leq \frac{1}{2} \\ 1 \text{ otherwise} \end{cases} \label{intro_bayes} \eeq
	is the Bayesian classifier under the distribution $Q$. The expectation $\E_Z$ in \eqref{intro_rate} is taken over the random realizations of all the observed data, namely the set $Z$, defined as
	\beq Z \triangleq \{(X_1^P, Y_1^P), ..., (X_{n_P}^P, Y_{n_P}^P), (X_1^Q, Y_1^Q), ..., (X_{n_Q}^Q, Y_{n_Q}^Q) \}. \label{intro_defZ} \eeq 
	
	Note that if one only had observations from the target distribution $Q$, under the same regularity conditions,  the minimax rate on the excess risk would be $n_Q^{-\frac{\beta(1+\alpha)}{2\beta+d}}$. Therefore, the additional term $n_P^{\frac{2\beta + d}{2\gamma\beta+d}}$ in the minimax rate \eqref{intro_rate} quantifies an ``effective sample size" for transfer learning from the source distribution $P$ relative to the target distribution $Q$, and $\frac{2\beta + d}{2\gamma\beta+d}$ can be viewed as an optimal rate of transfer. This result answers one of the main questions in transfer learning: $n_P^{\frac{2\beta + d}{2\gamma\beta+d}}$ is the total amount of information which can be transferred from $P$ to $Q$, and this quantity depends on the relative signal exponent $\gamma$ which characterizes the discrepancy between $P$ and $Q$ in posterior drift.
	
	We construct a two-sample weighted $K$-nearest neighbors ($K$-NN) classifier and show that it attains the optimal rate given in \eqref{intro_rate}. However, this classifier depends on the parameters $\alpha$, $\beta$, and $\gamma$, which are typically unknown in practice. In this paper, we also propose a data-driven classifier $\hat f_a$ that automatically adapts to the unknown model parameters $\alpha, \beta$ and $\gamma$, with an additional $\log$ term on the excess risk bound:
	\[ \sup_{(P,Q)\in\Pi} \E_Z \mathcal{E}_Q(\hat f_a) \lesssim  \p{\p{\frac{n_P}{\log (n_P + n_Q)}}^{\frac{2\beta + d}{2\gamma\beta+d}} + \frac{n_Q}{\log (n_P + n_Q)}}^{-\frac{\beta(1+\alpha)}{2\beta+d}}. \]
	This adaptive procedure is essentially different from either the non-adaptive procedure given in this paper, or any nonparametric classification procedures in the literature. The adaptive classifier is constructed based on the ideas inspired by Lepski's method for nonparametric regression. The construction begins with a small number of the nearest neighbors, and gradually increases the number of the neighbors used to make the decision. The algorithm terminates once an empirical  signal-to-noise ratio reaches a delicately designed threshold. It is shown that the resulting data-driven classifier automatically adapts to a wide collection of parameter spaces.
	
	In some applications, there are data available from multiple source distributions. Intuitively, the samples from all source distributions are helpful to the classification task under the target distribution. We also consider transfer learning in this setting under the posterior drift model. Suppose there are multiple source distributions $P_1,  \ldots, P_m$ and one target distribution $Q$, each pair of distributions $(P_i, Q)$ has a relative signal exponent $\gamma_i$, $i \in \{1, \ldots, m\}$. The minimax optimal rate of convergence is established and the result quantifies precisely the contributions from the data generated by the individual source distributions. An adaptive procedure is constructed and shown to simultaneously attain the optimal rate up to a logarithmic factor over a large class of parameter spaces.
	
	The rest of the paper is organized as follows. In Section \ref{Model.sec}, after some basic notations and definitions are introduced, the model for  transfer learning under the posterior drift model  is proposed in a nonparametric classification setting. In Section \ref{Minimax.sec}, we establish the minimax optimal rate by constructing a minimax optimal procedure with guaranteed upper bound and a matching lower bound. In section \ref{Adaptive.sec}, a data-driven adaptive classifier is proposed and is shown to adaptively attain the optimal rate of convergence, up to a logarithmic factor. Section \ref{Numerical.sec} investigates the numerical performance of the data driven procedure. In section \ref{application.sec}, a real data application is carried out to further illustrate the benifit of our method. Section \ref{Extension.sec} considers transfer learning with multiple source distributions and a brief discussion is given in Section \ref{Discussion.sec}. For reasons of space, we prove one main result in Section \ref{Proof.sec} and provide the proofs of the other results and some technical lemmas in the Supplementary Material \citep{CaiWei2019Supplement}.

\section{Problem Formulation}
\label{Model.sec}

In this section, we introduce the posterior drift model for transfer learning. We begin with notation and basic definitions.

	\subsection{Notations and definitions}
	
		For a distribution $G$, denote by $G(\cdot)$ and $\E_G(\cdot)$  respectively the probability and expectation under $G$. Denote by $P_X$ and $Q_X$ the marginal distribution of $X$ under the  joint distributions $P$ and $Q$ for $(X, Y)$ respectively. Let $\supp(\cdot)$ denote the support of a probability distribution. Throughout the paper we write $\Vert \cdot \Vert$ to denote the euclidean norm. We use $\ind{ \cdot }$ to denote the indicator function taking values in $\{0,1\}$. We define $a\vee b = \max(a,b)$, $a\wedge b = \min(a,b)$, and $\lfloor a \rfloor$ be the maximum integer that is not larger than $a$. We write $\lambda(\cdot)$ to denote Lebesgue measure of a set in a Euclidean space. We denote by $C$ or $c$ some generic constants not depending on $n_P$ or $n_Q$ that may vary from place to place.
	 
	\subsection{Posterior drift in nonparametric classification}
		
	For two distributions $P$ and $Q$ for a random pair $(X, Y)$ taking values in $[0,1]^d \times \{0, 1\}$, we observe two independent random samples, $(X_1^P, Y_1^P), \ldots, (X_{n_P}^P, Y_{n_P}^P) \iid P$ and $(X_1^Q, Y_1^Q), \ldots, (X_{n_Q}^Q, Y_{n_Q}^Q) \iid Q$. We shall use $P$-data and $Q$-data  to refer to the data sets drawn from the distributions $P$ and $Q$ respectively.  We consider the transfer learning problem when there is a posterior drift between $P$ and $Q$, i.e. the covariates/features $X$ are drawn from almost the same distributions (having the same support with bounded densities), but the response/label $Y$ has different conditional distributions given $X$ between $P$ and $Q$.
	
	The regression functions have been defined informally in the introduction, now we give a precise definition. Let 
	\[ \eta_P(x) = \begin{cases} P(Y=1|X=x) &\text{ if } x\in \supp(P_X) \\ \frac{1}{2} &\text{ otherwise} \end{cases} \]
	\[ \eta_Q(x) = \begin{cases} Q(Y=1|X=x) &\text{ if } x\in \supp(Q_X) \\ \frac{1}{2} &\text{ otherwise} \end{cases} \]
	denote the corresponding regression functions of $P$ and $Q$. Besides the previous definition $\eqref{intro_bayes}$ of Bayes classifier under the target distribution $Q$, We can similarly define the Bayes classifier for the source distribution $P$ as:
	\[ f^\ast_P(x) = \begin{cases} 0 \text{ if } \eta_P(x) \leq \frac{1}{2} \\ 1 \text{ otherwise} \end{cases}  .\]
	
	Now assume $(X^P,Y^P)$ is a data pair drawn from the distribution $P$. From the definition, given $X^P=x$, $Y^P$ is more likely to be equal to 1 if $f^\ast_P(x) = 1$ whereas $Y^P$ is more likely to be equal to 0 if $f^\ast_P(x) = 0$. It is similar for the distribution $Q$. Thus informally one can regard $f^\ast_P(x)$ ($f^\ast_Q(x)$) as the true label at the covariate value $x$ under the distribution $P$ ($Q$).
	
	In transfer learning, although the observed data are drawn from two or more different distributions, these distributions are usually related to each other so that all of them are useful for learning the intrinsic true labels. For instance, in a crowdsourcing survey, although accuracy varies among different workers, their answers should be no worse than random guessing. It is reasonable to assume that the answer is correct with probability at least $\frac{1}{2}$. This means we may reasonably assume that, given the same covariate $x$, the ``true labels" under the distributions $P$ and $Q$  are the same. That is	
	\[ f^\ast(x) \triangleq f^\ast_P(x) = f^\ast_Q(x) \quad \forall x\in \supp(P_X), \]
	which is equivalent to
	\[ (\eta_P(x) - \frac{1}{2})(\eta_Q(x) - \frac{1}{2}) \geq 0.  \]
	
	The definitions and assumptions introduced so far treat the $P$-data and $Q$-data symmetrically and interchangeably. But in real applications, usually the two data sets are treated differently. We call  $P$ the source distribution and $Q$ the target distribution. The goal is to transfer the knowledge gained from the $P$-data together with the information contained in the $Q$-data for constructing an optimal classifier under the target distribution $Q$. 
	
	Intuitively it is clear that the amount of information that can be transferred from the $P$-data for  the inference under $Q$ depends on the similarity between the distributions $P$ and $Q$.  In this paper, we quantify the similarity by the \textbf{relative signal exponent} of $P$ with respect to $Q$.
	
	\begin{definition}[Relative Signal Exponent]
		The class $\Gamma(\gamma, C_\gamma)$ with relative signal exponent $\gamma\in (0,\infty)$ and a constant $C_\gamma \in (0,\infty)$ is defined as the set of distribution pairs $(P,Q)$, both supported on $\reals^d \times \{ 0,1 \}$, satisfying $\forall x \in \supp(P_X)\cup \supp(Q_X)$, 	
		\beq (\eta_P(x) - \frac{1}{2})(\eta_Q(x) - \frac{1}{2}) \geq 0 \label{rse1} \eeq
		\beq |\eta_P(x) - \frac{1}{2}| \geq C_\gamma |\eta_Q(x) - \frac{1}{2}|^\gamma. \label{rse2} \eeq
	\end{definition}

\hspace{-10pt}
	\begin{figure}[t]
		\begin{minipage}{0.47\textwidth}
			\includegraphics[width=2.73in,height=2in]{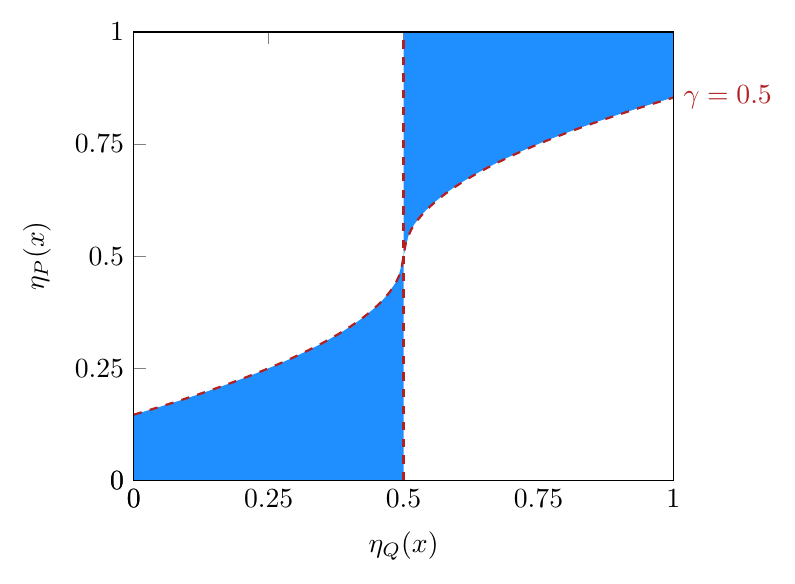}
		\end{minipage}
		\begin{minipage}{0.05\textwidth}
			\ 
		\end{minipage}
		\begin{minipage}{0.45\textwidth}
			\includegraphics[width=2.73in,height=2in]{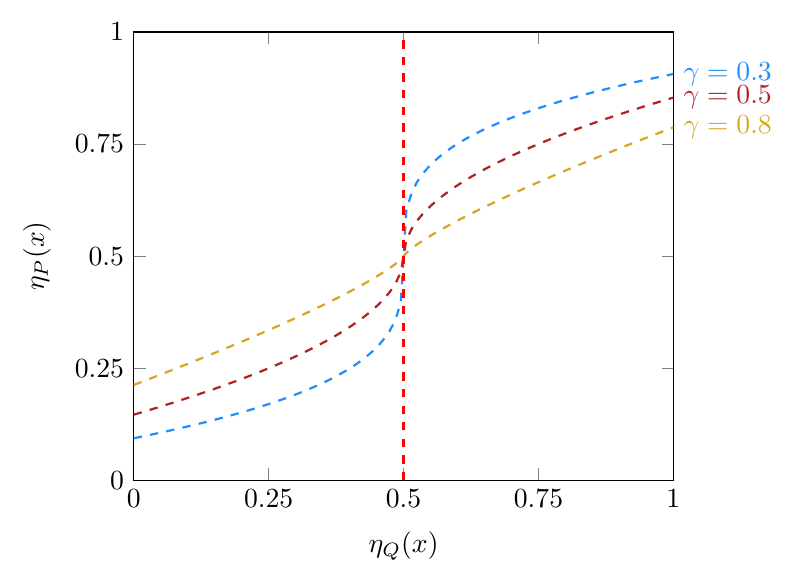}
		\end{minipage}
		\caption{\small Illustration of the relative signal exponent $\gamma$. Left panel: feasible region when $\gamma = 0.5$ and $C_\gamma = 0.5$. A pair of distributions $(P,Q)$ has relative signal exponent $\gamma = 0.5$ with $C_\gamma = 0.5$ when $(\eta_P(x), \eta_Q(x))$ falls into the shaded (blue) region for all $x$ in the support. Right panel: feasible region with different choices of $\gamma$. Smaller $\gamma$ implies more information contains in $P_{Y|X}$.}
	\end{figure}

	\begin{rmk}{ \rm
		The relative signal exponent $\gamma$ indicates the signal strength of the $P$-data relative to the $Q$-data. Note that $|\eta_Q(x) - \frac{1}{2}|$ is always bounded by $1/2$. So generally speaking, the smaller $\gamma$ is, the larger the difference between $\eta_P(x)$ and $\frac{1}{2}$, which means the $P$-data is more informative about $f^\ast(x)$ and consequently more information can be transferred from the $P$-data to the $Q$-data.
		
		One can see that the above definition of relative signal exponent implies when $|\eta_Q(x) - \frac{1}{2}|$ is large, then $|\eta_P(x) - \frac{1}{2}|$ should be relatively large. This is intuitively true in a wide range of real applications. Taking again the crowdsourcing surveys as an example. If one crowd of workers can answer a question correctly with a larger probability, then for another crowd of workers the accuracy of their answers is also usually larger because this question is likely to be easier.
}	\end{rmk}
	
	In addition to the relative signal exponent $\gamma$, we also need to define  a smoothness parameter of $\eta_Q$ and  characterize its behavior near $1/2$: 
	
	\begin{definition}[Smoothness]
		The $(\beta,C_\beta)-$H\"{o}lder class of functions ($0<\beta \leq 1$), denoted by $\mathcal{H}(\beta, C_\beta)$, is defined as the set of functions $g: \reals^{d} \to \reals$ satisfying, for any $x_1, x_2 \in \reals^d$,
		\[ |g(x_1) - g(x_2)| \leq C_\beta \Vert x_1 - x_2 \Vert^\beta. \]
	\end{definition}
	
	\begin{definition}[Margin Assumption]
		The margin class $\mathcal{M}(\alpha, C_\alpha)$ with $\alpha \geq 0, C_\alpha > 0$ is defined as the set of distributions $Q$ satisfying
		\[ Q_X(|\eta_Q(X) - \frac{1}{2}| < t) \leq C_\alpha t^\alpha. \]
	\end{definition}
	
	In this paper we consider the nonparametric classification problem when $\eta_Q(x)$ belongs to a $(\beta,C_\beta)-$H\"{o}lder class and $Q$ belongs to a margin class $\mathcal{M}(\alpha, C_\alpha)$. When $Q \in \mathcal{M}(\alpha, C_\alpha)$, we also say that $Q$ satisfies the margin assumption with the parameter $\alpha$.
	
	\begin{rmk}{\rm
		In the main part of our discussion, we focus on the case with $0<\beta\leq 1$, i.e. $\eta$ belongs to a H\"{o}lder function class with smoothness less than or equal to 1. Generally it is possible to consider more general classes where the smoothness parameter can be larger than 1. The discussion on the model and methods associated with the general smoothness parameter $\beta>1$ will be deferred to the discussion section.
		
		The margin assumption was first introduced in \cite{tsybakov2004optimal, audibert2007fast} to characterize the convergence rate in nonparametric classification. The margin assumption put a constraint on the mass around $\eta_Q(x) \approx \frac{1}{2}$ so that with large probability $\eta_Q(x)$ is either $\frac{1}{2}$ or far from $\frac{1}{2}$. Generally, if an underlying distribution satisfies the margin assumption, then a more accurate classification can be guaranteed.
}		
	\end{rmk}

	Another definition is about density constraints on the marginal distributions $P_X$ and $Q_X$.
	
	\begin{definition}[Common Support and Strong Density Assumption]
		$(P_X,Q_X)$ is said to have common support and satisfy strong density assumption with parameter $\mu = (\mu_-,\mu_+), c_\mu>0, r_\mu>0$ if both $P_X$ and $Q_X$ are absolutely continuous with respect to the Lebesgue measure on $\reals^d$, and
		\[ \Omega \triangleq \supp(P_X) = \supp(Q_X) \]
		\[ \lambda[\Omega \cap B(x,r)] \geq c_\mu\lambda[B(x,r)]  \quad \forall 0<r\leq r_\mu, \forall x\in \Omega \]
		\[ \mu_- < \frac{dP_X}{d\lambda}(x) < \mu_+ \quad \forall x\in \Omega \]
		\[ \mu_- < \frac{dQ_X}{d\lambda}(x) < \mu_+ \quad \forall x\in \Omega. \]
		
		Define $\mathcal{S}(\mu, c_\mu, r_\mu)$ to be the set of the marginal densities pairs $(P_X, Q_X)$ that have common support and satisfy the strong density assumption with parameter $\mu, c_\mu, r_\mu$.
	\end{definition}

	\begin{rmk}{\rm
		The strong density assumption was first introduced in \cite{audibert2007fast}. In this paper we focus on the scenario that the marginal densities of $P_X$ and $Q_X$ have regular support and are bounded from below and above on the support. 
		
		Moreover, note that when $Q_X$ satisfies the strong density assumption, in the regime $\alpha\beta > d$, there is no distribution $Q$ such that the regression function $\eta_Q$ crosses $\frac{1}{2}$ in the interior of the support $Q_X$ \citep{audibert2007fast}. So this regime only contains the trivial cases for classification.  Therefore, we further assume $\alpha\beta\leq d$ in the  following discussion. 
}	\end{rmk}

	Given a classifier $\hf: \R^d \to \{ 0,1 \}$, the excess risk on $Q$ of the classifier $\hf$, defined in equation \eqref{def_exrisk} in the introduction, has a dual representation \citep{gyorfi1978rate}
	\beq \mathcal{E}_Q(\hf) = 2\E_{(X,Y)\sim Q}(|\eta_Q(X) - \frac{1}{2}|\ind{\hf(X) = f^\ast_Q(X)}). \label{dual_rep} \eeq	
	
	A major goal in transfer learning is to construct an empirical decision rule $\hat f$ incorporating both the $P$-data and $Q$-data, so that the excess risk on $Q$ is minimized. It is interesting to understand when the minimax rate in the transfer learning setting is faster than the optimal rate  where only the $Q$-data is used to construct the decision rule.
	
	Putting the above definitions together, in this paper we consider the posterior drift nonparametric parameter space:
	\baligns 
	\Pi(\alpha, C_\alpha,  \beta,  C_\beta, \gamma, C_\gamma, \mu, c_\mu, r_\mu) = \{ (P,& Q) :  (P, Q)\in \Gamma(\gamma, C_\gamma), Q\in \mathcal{M}(\alpha, C_\alpha), \\
	&  \eta_Q\in\mathcal{H}(\beta,C_\beta),   (P_X,Q_X) \in \mathcal{S}(\mu, c_\mu, r_\mu)\}. 
	\ealigns	
	
	In the rest of this paper, we will use the shorthand $\Pi(\alpha,\beta,\gamma,\mu)$ or $\Pi$ if there is no confusion. The space $\Pi(\alpha,\beta,\gamma,\mu)$ is also called the posterior drift regime with $(\alpha,\beta,\gamma,\mu)$.

\section{Minimax Rate of Convergence}
\label{Minimax.sec}
	
	In this section, we establish the minimax rate of convergence for the excess risk on $Q$ for transfer learning under the posterior drift model and propose an optimal procedure using the two-sample weighted $K$-NN classifier. 
	
	The $K$-NN method has attracted much attention \citep{cover1967nearest, gyorfi1978rate, gadat2016} due to its massive practical success and appealing theoretical properties. In the conventional setting where one only has access to the $Q$-data and there is no $P$-data, with a suitable choice of the neighborhood size $k$,  the $K$-NN classifier can achieve the minimax rate of convergence for the excess risk on $Q$ \citep{gadat2016}. The $K-$NN classifier is generated in two steps:
	
	\begin{enumerate}
		\item[Step 1:] For any given $x$ to be classified, one can estimate $\eta_Q(x)$ by taking empirical mean of the response variables ($Y$) according to its $k$ nearest covariates ($X$). Formally, define $X_{(i)}^Q(x)$ be the $i$-th nearest covariates to $x$ among $X_1^Q,...,X_{n_Q}^Q$ and $Y^Q_{(i)}(x)$ is its corresponding response (label). The estimate $\hat\eta_Q(x)$ is given by
		\[ \hat\eta_Q(x) = \frac{1}{k} \sum_{i=1}^k Y_{(i)}^Q(x). \]
		
		\item[Step 2:] The class label for $x$ is estimated by the plug-in rule:
		\[  \hat f(x) = \ind{\hat\eta_Q(x) > \frac{1}{2}}. \]
	\end{enumerate}
	
	In transfer learning, one also has access to the $P$-data in addition to the $Q$-data, the $P$-data can be used to help the classification task under the target distribution $Q$ and should be taken into consideration. To accommodate the existing $K$-NN methods, we should take the empirical mean of not only the $k$-nearest response variables from the $Q$-data, but also some nearest response variables from the $P$-data. In addition, when taking the average, data from the different distributions should have different weights because the signal strength varies between the two distributions. To make the classification at $x \in [0,1]^d$, a new strategy called the two-sample weighted $K$-NN classifier is summarized as follows:
	
	\begin{enumerate}
		\item[Step 1:] Define $X_{(i)}^P(x)$ to be the $i$-th nearest covariates to $x$ among $X_1^P,...,X_{n_P}^P$ and $Y^P_{(i)}(x)$ is its corresponding response. $X_{(i)}^Q(x)$ and $Y^Q_{(i)}(x)$ can be defined likewise. Construct the two-sample weighted $K$-NN estimator
		\[ \hat{\eta}_{NN}(x) = \frac{w_P\sum_{i=1}^{k_P}Y^P_{(i)}(x) + w_Q\sum_{i=1}^{k_Q}Y^Q_{(i)}(x)}{w_Pk_P + w_Qk_Q} \]
		
		where the number of neighbors $k_P$ and $k_Q$ and the weights $w_P$ and $w_Q$ will be specified later.
		
		\item[Step 2:] The class label for $x$ is estimated by the plug-in rule:
		\[ \hf_{NN}(x) = \ind{ \hat\eta_{NN}(x)>\frac{1}{2} }. \]
	\end{enumerate}
 
	The final decision rule $\hf_{NN}(x)$ is called the {\it two-sample weighted $K$-NN classifier.} It is clear that  $\hf_{NN}(x)$ is generated by both the $P$-data and $Q$-data.
	
	The performance of the  two-sample weighted $K$-NN classifier $\hf_{NN}(x)$ clearly depends on the choice of $(k_P, k_Q, w_P, w_Q)$. The next theorem gives a set of choices of $(k_P, k_Q, w_P, w_Q)$ and a provable upper bound on the excess risk, which gives a guarantee for the performance of the two-sample weighted $K$-NN classifier with these specific choices of $(k_P, k_Q, w_P, w_Q)$.
	
	\begin{theorem}[Upper Bound] \label{ub}
		 Let $\hat f_{NN}$ be the  two-sample weighted $K$-NN classifier with $w_Q = (n_P^{\frac{2\beta+d}{2\gamma\beta+d}}+n_Q)^{-\frac{\beta}{2\beta+d}}$, $w_P = (n_P^{\frac{2\beta+d}{2\gamma\beta+d}}+n_Q)^{-\frac{\gamma\beta}{2\beta+d}}$, $k_Q = \lfloor n_Q(n_P^{\frac{2\beta+d}{2\gamma\beta+d}}+n_Q)^{-\frac{d}{2\beta+d}} \rfloor$, and $k_P = \lfloor n_P(n_P^{\frac{2\beta+d}{2\gamma\beta+d}}+n_Q)^{-\frac{d}{2\beta+d}}\rfloor$. Then 
		\[ \sup_{(P,Q) \in \Pi} \E \mathcal{E}_Q(\hf_{NN}) \leq C(n_P^{\frac{2\beta + d}{2\gamma\beta+d}} + n_Q)^{-\frac{\beta(1+\alpha)}{2\beta+d}} \]
		for some constant $C > 0$ not depending on $n_P$ or $n_Q$.
	\end{theorem}

	The following lower bound result shows that the two-sample weighted $K$-NN classifier $\hf_{NN}$ given in Theorem \ref{ub} is in fact rate optimal.
	
	\begin{theorem}[Lower Bound] \label{lb}
	There exists a constant $c > 0$ not depending on $n_P$ or $n_Q$ such that
	\[ \inf_{\hat f}\sup_{(P,Q) \in \Pi} \E \mathcal{E}_Q(\hat{f}) \geq c(n_P^{\frac{2\beta + d}{2\gamma\beta+d}} + n_Q)^{-\frac{\beta(1+\alpha)}{2\beta+d}}. \]
	\end{theorem}
	
	\begin{figure}[t]
		\centering
		\includegraphics[width=0.96\textwidth]{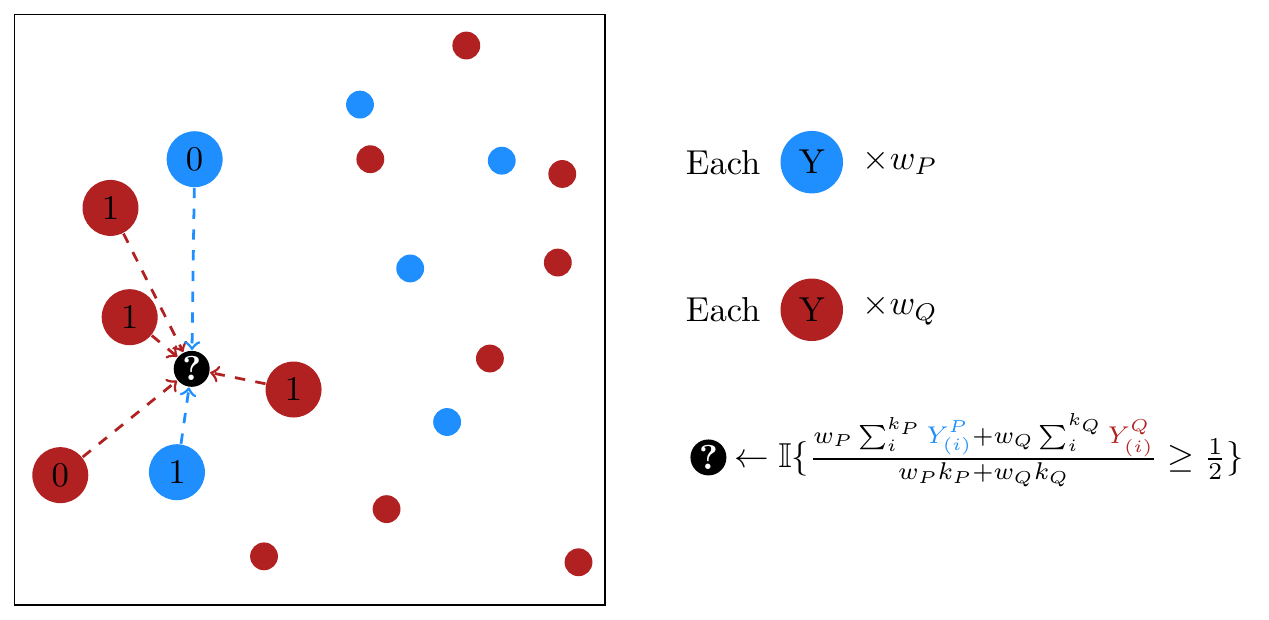}
		\caption{\small An illustration of the two-sample weighted $K$-NN classifier. $(X^P,Y^P)$ are shown by the blue points and $(X^Q,Y^Q)$ are shown by the red points. For each point in the graph, the coordinates represent its two-dimensional covariates $X$ while the number marked inside the point represents its label $Y$. To classify the black point ($x$) located in middle of the graph, by calculation we get (say) $k_P=2$ and $k_Q=4$. Then we find $k_P$ nearest neighbors from $P$-data and $k_Q$ nearest neighbors from $Q$-data. Finally, we calculate their weighted mean to make the final classification. }
		\label{nonada_pic}
	\end{figure}
	
	The proof of Theorem \ref{ub} will be given in Section \ref{Proof.sec}. The proof of Theorem \ref{lb} will be given in the supplementary material \citep{CaiWei2019Supplement}. Theorems \ref{ub} and \ref{lb} together establish the minimax rate of convergence for transfer learning under the posterior drift model,
	\beq \inf_{\hat f}\sup_{(P,Q) \in \Pi} \E \mathcal{E}_Q(\hat{f}) \asymp (n_P^{\frac{2\beta + d}{2\gamma\beta+d}} + n_Q)^{-\frac{\beta(1+\alpha)}{2\beta+d}}. \label{rate} \eeq
	
	We make a few remarks on the minimax rate of convergence.
	\begin{itemize}
	\item Based on the minimax  rate given in \eqref{rate}, it is easy to see that, in terms of the classification accuracy,  the contribution from the $P$-data is substantial when $n_P^{\frac{2\beta + d}{2\gamma\beta+d}} \gg n_Q$, and the contribution is not significant otherwise.

	\item It is worth noting that in the conventional setting with access to the $Q$-data only,  the minimax rate, which is given in \cite{audibert2007fast}, would be
		\beq  \inf_{\hat f}\sup_{(P,Q) \in \Pi} \E \mathcal{E}_Q(\hat{f}) \asymp n_Q^{-\frac{\beta(1+\alpha)}{2\beta+d}}, \label{rateq}  \eeq
which is a special case of \eqref{rate} with $n_P=0$. This rate can be achieved by the $K$-NN classifier given above with the choice of $k \asymp n_Q^{\frac{2\beta}{2\beta+d}}$.
				
	\item Comparing the convergence rates \eqref{rate} with \eqref{rateq}, the minimax rate for transfer learning under the posterior drift model is the same as if one had a sample of size $n_P^{\frac{2\beta + d}{2\gamma\beta+d}} + n_Q$ from the distribution $Q$ in the conventional setting. Therefore, one can intuitively view $n_P^{\frac{2\beta + d}{2\gamma\beta+d}}$ as the ``effective sample size" of the $P$-data for the  classification task under $Q$. The exponent $\frac{2\beta + d}{2\gamma\beta+d}$ here can be regarded as the {\it transfer rate}. The smaller the relative signal exponent $\gamma$ is, the larger $\frac{2\beta + d}{2\gamma\beta+d}$ is, and more information is transferred from the $P$-data. This transfer rate provides a quantitative answer to the question posed  in the introduction: How much information can be transferred from the source distribution $P$ to the target distribution $Q$? It is  also interesting to note that, when $\gamma<1$,  $\frac{2\beta + d}{2\gamma\beta+d}>1$, which implies that  in this case an observation from $P$ is more valuable than an observation from $Q$ for the classification problem. 
		
	\item	In the transfer learning literature, much attention has been on an interesting special case where there is no data from the target distribution $Q$ at all, i.e., $n_Q = 0$ \citep{mansour2009domain, blitzer2008learning}. This case arises when a classifier has been trained based on the data drawn from the source distribution $P$, and one wishes to generalize the classifier to unlabeled testing data drawn from the target distribution $Q$. Our results show that generalization is possible in the posterior drift  framework and the optimal rate of convergence is
		\[ \inf_{\hat f}\sup_{(P,Q) \in \Pi} \E \mathcal{E}_Q(\hat{f}) \asymp n_P^{-\frac{\beta(1+\alpha)}{2\gamma\beta+d}}. \]
		
	
	\end{itemize}

\section{Data-driven Adaptive Classifier}
\label{Adaptive.sec}
	
	In the previous section, we have established the minimax optimal rate over the parameter space $\Pi(\alpha,\beta,\gamma,\mu)$ for transfer learning under the posterior drift model. This rate can be achieved by the two-sample weighted $K$-NN classifier given in Theorem \ref{ub}. A major drawback of this classifier is that it requires the prior knowledge of $\beta$ and $\gamma$, which is typically unavailable in practice. An interesting and practically important question is whether it is possible to construct a data-driven adaptive decision rules that can achieve the same rate of convergence, while automatically adapt to a wide collection of the parameter spaces $\Pi(\alpha,\beta,\gamma,\mu)$.
	
	In nonparametric regression, Lepski's method  \citep{lepskii1991problem,lepskii1992asymptotically,lepskii1993asymptotically}
is a well known approach for the construction of a data driven estimator that adapts to the unknown smoothness parameter $\beta$ by screening from a small bandwidth to larger bandwidths with delicately designed stopping rules. This procedure can be used for nonparametric classification in the conventional setting where only $Q$-data is available and only adaptation to one unknown smoothness parameter $\beta$ is needed. For readers' convenience we include this construction in Section \ref{Proof.sec}. The transfer learning setting is more challenging: we need to adapt to both unknown parameters $\beta$ and $\gamma$. In this section, we modify Lepski's method suitably in our context and introduce a new stopping rule and show that the resulting data driven classifier can adapt to all unknown parameters.

	Now we develop a data-driven procedure to make classification at a specific point $x\in [0,1]^d$. In our construction, we need to combine all data points from the $P$-data and the $Q$-data together and find nearest neighbors among all the data. Denote by $X_{(i)}(x)$ the $i$-th nearest data point to $x$ in the combined set $\{ X_1^P, ..., X_{n_P}^P \} \cup \{ X_1^Q, ..., X_{n_Q}^Q \}$.  Similar to Lepski's method, we begin with a small number of nearest neighbors, and gradually increase the number of neighbors used to make the decision. One more nearest neighbor is added  in each step. At the $k$-th step, there are $k$ nearest neighbors $X_{(1)}(x), ..., X_{(k)}(x)$ among all the points in the combined set $\{ X_1^P,  ..., X_{n_P}^P \} \cup \{ X_1^Q, ..., X_{n_Q}^Q \}$. Suppose among these $k$ nearest neighbors  there are $k_P^{(k)}$ points from the $P$-data and $k_Q^{(k)}$ points from the $Q$-data.
	Heuristically, given these $k$ nearest neighbors, one can obtain a weighted $k$-NN estimate as
	\[ \hat{\eta}^{(k)}(x, w_P, w_Q) = \frac{w_P\sum_{i=1}^{k_P^{(k)}}Y^P_{(i)}(x) + w_Q\sum_{i=1}^{k_Q^{(k)}}Y^Q_{(i)}(x)}{w_Pk_P^{(k)} + w_Qk_Q^{(k)}}. \]
	
	If $\beta$ and $\gamma$ are known, one can calculate the optimal choice of the weights $w_P$ and $w_Q$ for a two-sample weighted $K$-NN classifier. To construct an adaptive procedure, we need to find a data driven method for choosing the weights $w_P$ and $w_Q$. Define the ``variance" of $\hat{\eta}^{(k)}(x, w_P, w_Q)$ as
	\[ v^{(k)}(w_P, w_Q) = \frac{w_P^2k_P^{(k)}+w_Q^2k_Q^{(k)}}{(w_Pk_P^{(k)} + w_Qk_Q^{(k)})^2}.  \] 
	
	For a given $k$, we call the maximum value of the ratio between $(\hat{\eta}^{(k)}(x, w_P, w_Q)- \frac{1}{2})^2$ and the ``variance" $v^{(k)}(w_P, w_Q)$ as the signal-to-noise ratio index $\hat{r}^{(k)}$:
	\[ \hat{r}^{(k)} = \max_{w_P,w_Q} \frac{(\hat{\eta}^{(k)}(x, w_P, w_Q)- \frac{1}{2})^2}{ v^{(k)}(w_P, w_Q)}. \]
	
	The algorithm is terminated when $\hat{r}^{(k)} > (d+3)\log(n_P+n_Q)$, and the corresponding $w_p$ and $w_Q$ are chosen as the maximizers of $\frac{(\hat{\eta}^{(k)}(x, w_P, w_Q)- \frac{1}{2})^2}{v^{(k)}(w_P, w_Q)}$. If the algorithm doesn't terminate at any step, the optimal $k$ can be alternatively chosen by the maximizer of $\hat{r}^{(k)}$. That is, we choose $k=k^*$ with
	\beq k^* = \begin{cases} \min\{k:  \hat{r}^{(k)} > (d+3)\log(n_P+n_Q)\} &\text{ if } \max_k \hat{r}^{(k)} >  (d+3)\log(n_P+n_Q) \\ \argmax_k \hat{r}^{(k)} \quad &\text{ otherwise} \end{cases} \label{stopping-rule}  \eeq
and choose $ (w_P, w_Q) = (w_P^\ast, w_Q^\ast) $ with 
	\[ (w_P^\ast, w_Q^\ast) = \argmax_{(w_P,w_Q)} \frac{(\hat{\eta}^{(k^\ast)}(x, w_P, w_Q)- \frac{1}{2})^2}{v^{(k^\ast)}(w_P, w_Q)}. \]
The data driven adaptive classifier is then defined as
	\[ \hf_a(x) = \ind{\hat{\eta}^{(k^*)}(x, w_P^\ast, w_Q^\ast) \geq \frac{1}{2}}. \]

\begin{rmk}{\rm	
	The choice of $(d+3)\log(n_P+n_Q)$ as the threshold in the stopping rule \eqref{stopping-rule} is important and requires some explanation. Roughly speaking, this is due to the fact that the maximum fluctuation of $\hat{\eta}^{(k)}(x, w_P, w_Q)$ is bounded by $\sqrt{(d+3)\log(n_P+n_Q)v^{(k)}(w_P, w_Q)}$ with high probability, which will be shown in Lemma 5 with a suitable change of parameter (stated in the supplementary material \citep{CaiWei2019Supplement}). Thus, when $\hat{r}^{(k)} > (d+3)\log(n_P+n_Q)$, $\hat{\eta}^{(k)}(x, w_P, w_Q) > \frac{1}{2}$ indicates $\E\hat{\eta}^{(k)}(x, w_P, w_Q) > \frac{1}{2}$, which suggests $f^\ast(x) = 1$, and vice versa. 
}\end{rmk}	

	The procedure is summarized in Algorithm \ref{ada_algo}. In the algorithm below we simplify the above procedure by providing the actual closed form expression for $\hat r^{(k)}$ and $\hf_a(x)$. 
	
	\begin{algorithm}[t]
		\caption{The Data Driven Procedure}
		\begin{algorithmic}
			\State \textbf{Input: } $x \in \supp(Q_X)$
			\For{$k = 1, ...,(n_P+n_Q-1),(n_P+n_Q) $}
				\State Find $k$ nearest covariates to $x$ among all covariates in data $\{ X_1^P, X_2^P, ..., X_{n_P}^P \} \cup \{ X_1^Q, X_2^Q, ..., X_{n_Q}^Q \}$. Suppose among those $k$ nearest neighbors $X_{(1)}(x), X_{(2)}(x), ..., X_{(k)}(x)$ there are $k_P^{(k)}$ covariates from $P$-data and $k_Q^{(k)}$ covariates from $Q$-data.
				\State Compute $k_P^{(k)}$ nearest neighbor estimate in $P$-data (If $k_P^{(k)}=0$, set $\hat\eta_P^{(k)} \gets \frac{1}{2}$ )
				\[ \hat\eta_P^{(k)} \gets \frac{1}{k_P^{(k)}}\sum_{i=1}^{k_P^{(k)}} Y^P_{(i)}(x)\]
				
				and $k_Q^{(k)}$ nearest neighbor estimate in $P$-data (If $k_Q^{(k)}=0$, set $\hat\eta_Q^{(k)} \gets \frac{1}{2}$ )
				\[ \hat\eta_Q^{(k)} \gets \frac{1}{k_Q^{(k)}}\sum_{i=1}^{k_Q^{(k)}} Y^P_{(i)}(x)\]
				\State Let $\hat r^{(k)}$ be the signal-to-noise ratio index calculated by
				\[ \hat r^{(k)} \gets \begin{cases} k_P^{(k)}\p{\hat\eta_P^{(k)} - \frac{1}{2}}^2 + k_Q^{(k)} \p{\hat\eta_Q^{(k)}- \frac{1}{2}}^2 \text{ if } \sign(\hat\eta_P^{(k)} - \frac{1}{2}) = \sign(\hat\eta_Q^{(k)} - \frac{1}{2}) \\ \max\p{k_P^{(k)}\p{\hat\eta_P^{(k)}- \frac{1}{2}}^2, k_Q^{(k)} \p{\hat\eta_Q^{(k)}- \frac{1}{2}}^2} \text{ if } \sign(\hat\eta_P^{(k)}- \frac{1}{2}) \neq \sign(\hat\eta_Q^{(k)}- \frac{1}{2}) \end{cases}  \]
				\State Define the intermediate classifier  by
				\[  \hat f^{(k)}(x) \gets \ind{\sqrt{k_P^{(k)}}\p{\eta_P^{(k)} - \frac{1}{2}} + \sqrt{k_Q^{(k)}} \p{\eta_Q^{(k)}- \frac{1}{2}} \geq 0} \]
				\If {$\hat r^{(k)}(x) > (d+3)\log(n_P+n_Q)$}
					\State	Stop and output $\hat f_a(x) \gets f^{(k)}(x)$
				\EndIf
			\EndFor
			
			\State Output $\hat f_a(x) \gets \hf^{(k_m)}(x)$ where $k_m = \argmax_k \hat r^{(k)}$
		\end{algorithmic}
		\label{ada_algo}
	\end{algorithm}
	
	\begin{figure}[t]
		\centering
		\includegraphics[width=0.96\textwidth]{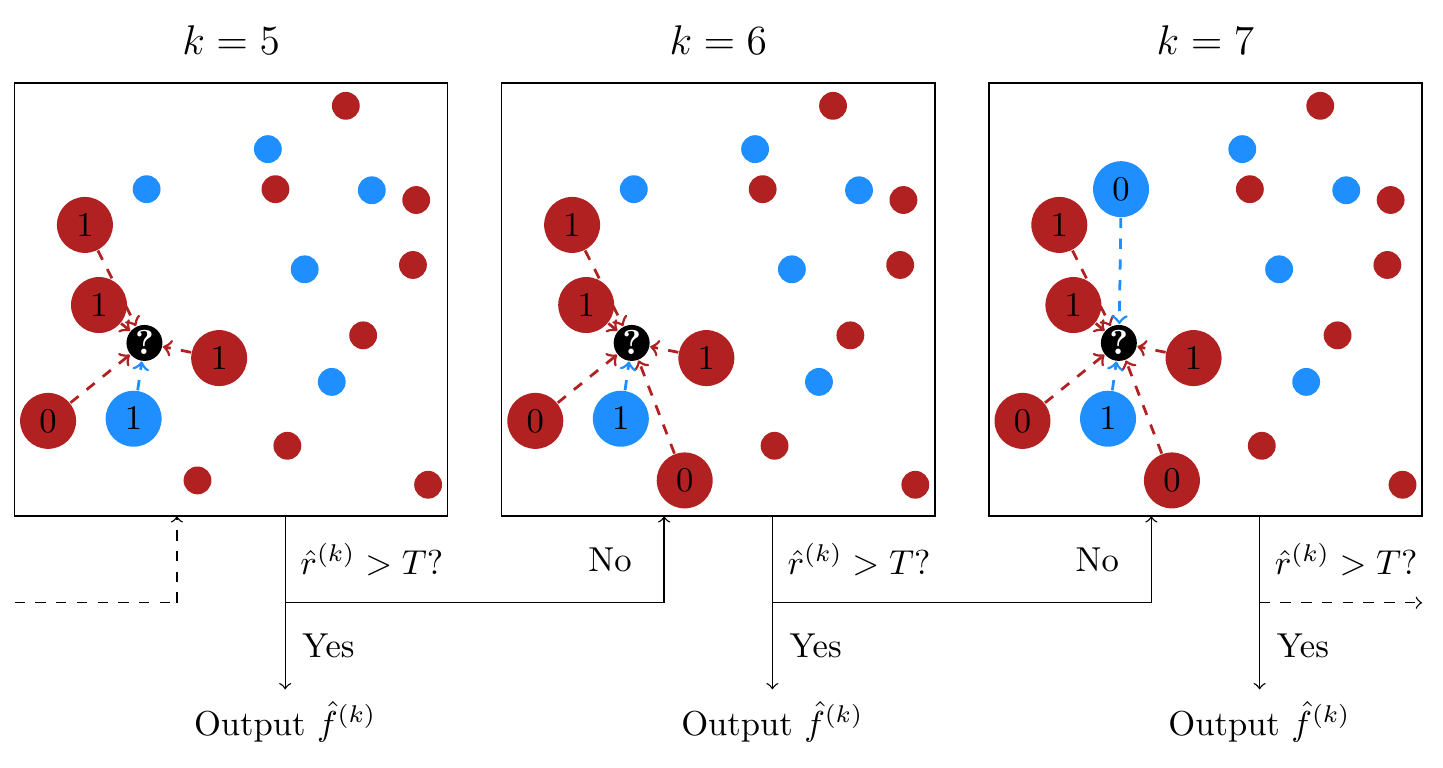}
		\caption{An illustration of the adaptive procedure given in Algorithm \ref{ada_algo}. See figure \ref{nonada_pic} for interpretation of the graph. Here we shorthand the threshold $T = (d+3)\log(n_P+n_Q)$. In each step, we evaluate $r^{(k)}$ and compare it to the threshold $R$. If $r^{(k)} > T$, then output $\hat f^{(k)}$ generated in current step; if $r^{(k)} \leq T$, go to next step and add one more nearest neighbor.}
		\label{ada_pic}
	\end{figure}
	
	Note that $\hf_a$ is a data-driven adaptive decision rule. For this adaptive classifier, the following theorem gives an upper bound for the excess risk under the target distribution $Q$:
	
	\begin{theorem} \label{adap_bound}
		Let $n=n_P+n_Q$. There exists a constant $C>0$ not depending on $n_P$ or $n_Q$ such that 
		\beq \sup_{(P,Q) \in \Pi} \E \mathcal{E}_Q(\hat{f}_a) \leq C \p{\p{\frac{n_P}{\log n}}^{\frac{2\beta + d}{2\gamma\beta+d}} + \frac{n_Q}{\log n}}^{-\frac{\beta(1+\alpha)}{2\beta+d}}. \label{adaub} \eeq
	\end{theorem}
	The proof of  Theorem \ref{adap_bound} is given in the supplementary material \citep{CaiWei2019Supplement}.
	
	Comparing  the rate of convergence \eqref{adaub} for the adaptive classifier $\hat f_a$ with the minimax rate \eqref{rate}, the data driven classifier $\hat f_a$ simultaneously achieves within a logarithmic factor of the minimax optimal rate for a large collection of parameter spaces.  The additional logarithmic term appearing in the rate is a common phenomenon in adaptive nonparametric regression under the pointwise loss. See, for example, \cite{lepskii1991problem} and \cite{brown1996constrained}.
	
	\begin{rmk}{\rm
		If we apply Lepski's method in the conventional setting where only the $Q$-data is available, then we have the following upper bound on the excess risk under $Q$:
		\beq \sup_{(P,Q) \in \Pi} \E \mathcal{E}_Q(\hat{f}_L) \leq C \cdot  \p{ \frac{n_Q}{\log  n_Q}}^{-\frac{\beta(1+\alpha)}{2\beta+d}}.\label{lepskib} \eeq	
		One can verify that by setting $n_P = 0$, our new adaptive procedure is exactly equivalent to Lepski's method, while the rates of convergence for the two methods also coincide. 
}	\end{rmk}

\section{Numerical Studies}
\label{Numerical.sec}

	In this section, we carry out simulation studies to further illustrate the performance of the adaptive transfer learning procedure.  Numerical comparisons with the existing methods are given. The simulation results are consistent with the theoretical predictions.
	
	For all simulation experiments in this section, we generate our data under the posterior drift model where $d=2$, $\eta_Q\in \mathcal{H}(1, 1)$,  $(P,Q)\in \Gamma(0.3, 1)$ with the relative signal exponent $\gamma = 0.3$, $Q$ satisfies the margin assumption with $\alpha = 0$, and  $P_X$ and  $Q_X$ have the common support and bounded densities.

	Now we specify the distribution $(P,Q)$ used to generate data. Let $x_c = (0.5,0.5)$. We set $P_X = Q_X$ to be both uniformly distributed on $[0,1]^3$. In addition, we set $\eta_Q$ and $\eta_P$ as
	\[ \eta_Q(x) = \max (p_{\max} - \Vert x-x_c \Vert, 0.5)  \quad{\rm and} \quad \eta_P(x) = 0.5+ (\eta_Q(x) - 0.5)^{0.3} \]	
	where $p_{\max} \in (0.5,1]$ is a variable in our simulation. According to the construction, $\eta_P$ and $\eta_Q$ both take maximum value at $x_c$ and are decreasing as $x$ moving away from $x_c$, then remain at a constant value 0.5 when $\Vert x-x_c \Vert \geq p_{\max} - 0.5$. Note that $\eta_Q(x_c)=p_{\max}$, so $p_{\max}$ indicates the difficulty to make classification on $x_c$, and we will use $p_{\max}$ and $\eta_Q(x_c)$ interchangeably in the following discussion. One can verify that the above construction of $(P,Q)$ satisfies the conditions we gave.
	 
	In the following experiments, we focus on evaluating the average classification accuracy at the random test sample $x$ drawn from the uniform distribution on the ball $B(x_c, 0.05)$, given $n_P$ data generated from the above distribution $P$ and $n_Q$ data generated from $Q$. This is arguably a good alternative to drawing the sample point from $Q_X$, because we always set $0.5 < p_{\max} \leq 0.55$ during the following experiments so that $\eta(x) = 0.5$ for all $x \notin B(x_c, 0.05)$. Compared to drawing from $Q_X$, the accuracy is scaled larger in order to better illustrate the results.

	\subsection{Minimax non-adaptive classifier}
	
	For this particular distribution pair $(P,Q)$, the minimax rate of convergence for the  excess risk can be achieved via two-sample weighted $K$-NN classifier, with choice of parameters $\gamma = 0.3$ and $\beta = 1$. We compare the accuracy of our classifier compared with two alternative methods. 
	
	The first alternative method is using a $K$-NN classifier on the combined dataset, ignoring the fact that they are from the different distributions. This method is the most commonly used  method in practice when people do not know how to make best transfer of information from $P$ into $Q$. 	
	The second alternative method is using a $K$-NN classifier on just the $Q$-data. By comparing this method with the proposed method one can know how much can be gained from the transfer of information.
	
\begin{figure}[t]
	\begin{subfigure}[b]{0.495\textwidth}
		\includegraphics[width=2.25in,height=1.84in]{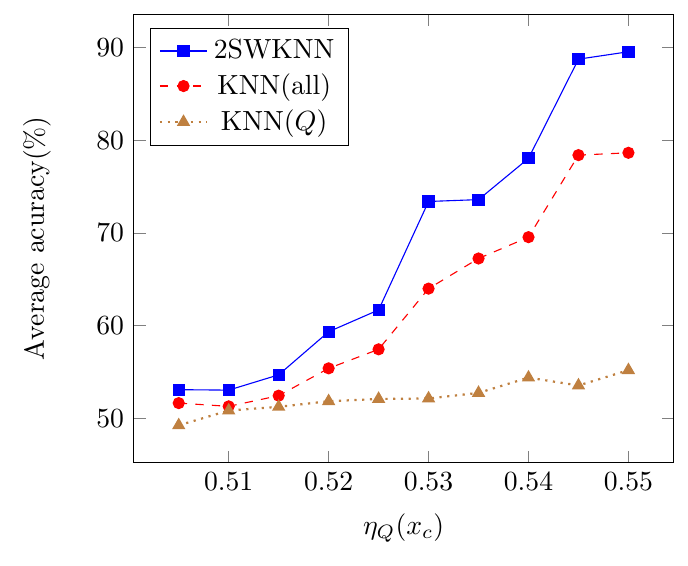}
	\caption{\small Experiments on different $\eta_Q(x_c)$}
	\label{nonadap_simu_a}
	\end{subfigure}
	\begin{subfigure}[b]{0.495\textwidth}
		\includegraphics[width=2.25in,height=1.84in]{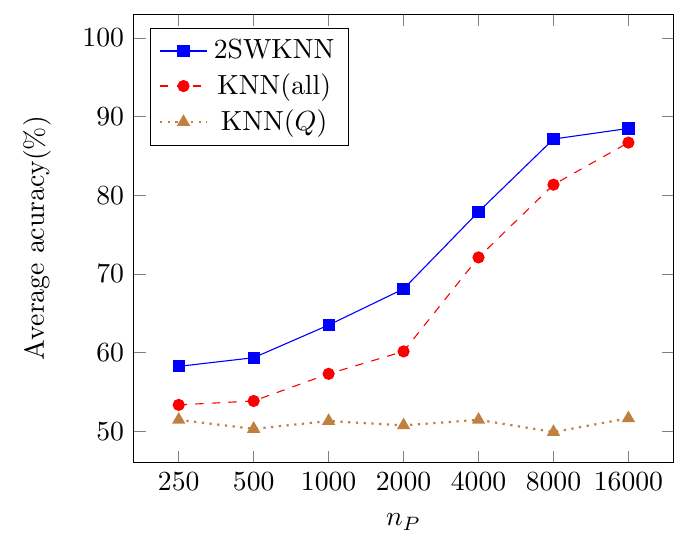}
	\caption{\small Experiments on different $n_p$}
	\label{nonadap_simu_b}
	\end{subfigure}
	\caption{\small Accuracy of non-adaptive classifiers. Blue: Two-sample weighted $K$-NN classifier. Red: $K$-NN using combined data. Brown: $K$-NN using only data from distribution $Q$. } \vspace{-15pt}
	\end{figure}

	With $n_P = 2000$ and $n_Q = 5000$,  Figure \ref{nonadap_simu_a} shows the average accuracy of classification on random test sample $x$ versus different $\eta_Q(x_c)$ during 2000 rounds of simulations, where $\eta_Q(x_c) \in \{0.505, 0.510, 0.515 ..., 0.55\}$. The further $\eta_Q(x_c)$ is away from $1/2$, the easier this classification task becomes. It can be seen from the plot that there is an obvious gap between the performance of our method and other methods, as the optimality theory predicts.
	
	With $n_Q = 5000$ and $\eta_Q(x_c) = 0.53$, Figure \ref{nonadap_simu_b} plots the accuracy versus different $n_P$ with 1000 replications, where $n_P \in \{250, 500, ..., 16000\}$. It can be seen from the plot that as  $n_P$ increases, the accuracy of the proposed method also improves because more information is available from the $P$-data to help  the classification task under the distribution $Q$. Also, a clear gap is seen on most values of $n_P$ between our classifier and other methods.
	
	\subsection{Adaptive classifier}
	\label{numerical.ada.sec}
	
	We also compare the proposed adaptive classifier with the existing methods to see whether its numerical  performance matches its theoretical guarantees. Lepski's method is a good competitor as it is also adaptive to the smoothness parameter $\beta$. Following a similar routine as in the previous experiments, we compare the classification accuracy for 3 methods: our proposed classifier, Lepski's method with all data involved, and Lepski's method with only the $Q$-data.
	
\begin{figure}[t]
		\begin{subfigure}[b]{0.495\textwidth}
			\includegraphics[width=2.25in,height=1.84in]{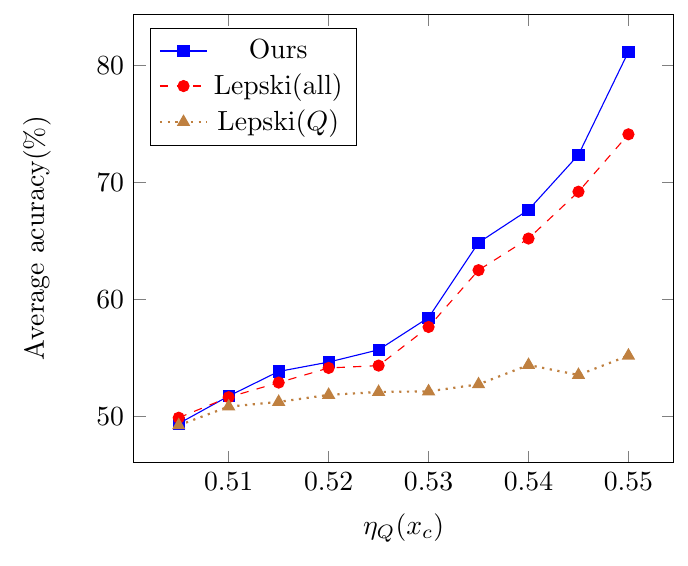}
			\caption{\small Experiments on different $\eta_Q(x_c)$}
			\label{adap_simu_a}
		\end{subfigure}
		\begin{subfigure}[b]{0.495\textwidth}
			\includegraphics[width=2.25in,height=1.84in]{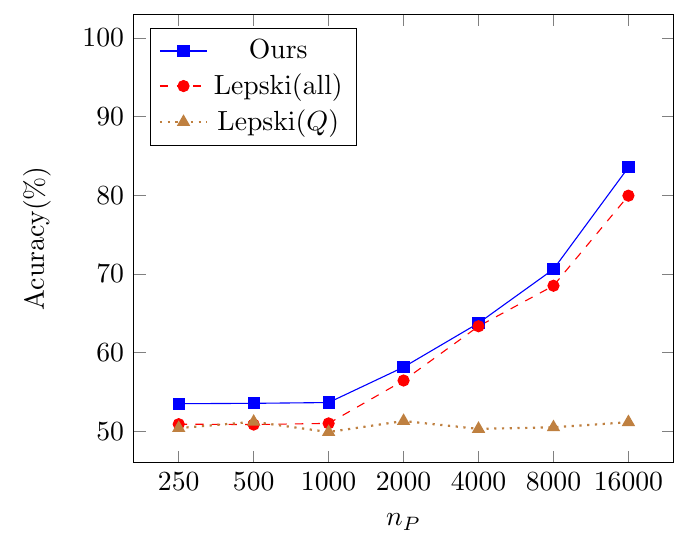}
			\caption{\small Experiments on different $n_p$}
			\label{adap_simu_b}
		\end{subfigure}
		\caption{\small Accuracy of adaptive methods. Blue: Our adaptive classifier. Red: Lepski's method using combined data. Brown: Lepski's method using only data from distribution $Q$. }  \vspace{-15pt}
	\end{figure}

	Figure \ref{adap_simu_a} has exactly the same setup as in Figure \ref{nonadap_simu_a} with $n_P = 2000$ and $n_Q = 5000$. We plot the classification accuracy  at $x_c$ versus different $\eta_Q(x_c)$ with 2000 replications, where $\eta_Q(x_c) \in \{0.505, 0.51, 0.515 ..., 0.55\}$. Although the difference is smaller, there is still a gap between our method and Lepski's methods. The gap increases as $n_Q(x_c)$ increases. Compared with the non-adaptive methods, the proposed adaptive classifier does not lose much accuracy without the knowledge of $\gamma$ and  $\beta$.
	 
	With $n_Q = 5000$ and $\eta_Q(x_c) = 0.53$, Figure \ref{adap_simu_b} plots the accuracy versus different $n_P$ with 1000 replications, where $n_P \in \{250, 500, ..., 16000\}$. Our adaptive classifier consistently outperforms the other methods for all $n_P$.
	
\section{Application to Crowdsourced Mapping Data}
\label{application.sec}

To illustrate the proposed adaptive classifier, we consider in this section an application based on the crowdsourced mapping data \citep{johnson2016integrating}. Land use/land cover maps derived from remotely-sensed imagery are important for geographic studies. This dataset contains Landsat time-series satellite imagery information on given pixels and their coresponding land cover class labels (farm, forest, water, etc.) obtained from multiple sources. The goal is to make classification of land cover classes based on NDVI (normalized difference vegetation index) values of those remotely-sensed imagery from the years 2014-2015. In this paper we focus on classification of two specific classes: farm and forest. 

Within this dataset, there are two kinds of label sources, given the NVDI values of the images: 1) crowdsourced georeferenced polygons with land cover labels obtained from OpenStreetMap; 2) accurately labeled data by experts. Although crowdsourced data are massive, free and public, the labels contain various types of errors due to user mislabels or outdated images. Whereas the expert labels are almost accurate, but they are usually too expensive to obtain a large volume. The challenge is to accurately combine the information contained in the two datasets to minimize the classification error. 

As in Section \ref{numerical.ada.sec}, we apply three methods to make the classification: the proposed adaptive procedure, Lepski's method with all data involved, and Lepski's method with only the crowdsourced data. We use $n_P = 50$ accurately labeled data, and change the number of involved crowdsourced data from $n_Q = 25$ to $n_Q = 800$. We use other 166 accurately labeled data to evaluate the classification accuracy of the three methods mentioned above.

\begin{figure}[t]
	\begin{subfigure}[b]{0.495\textwidth}
		\includegraphics[width=2.6in]{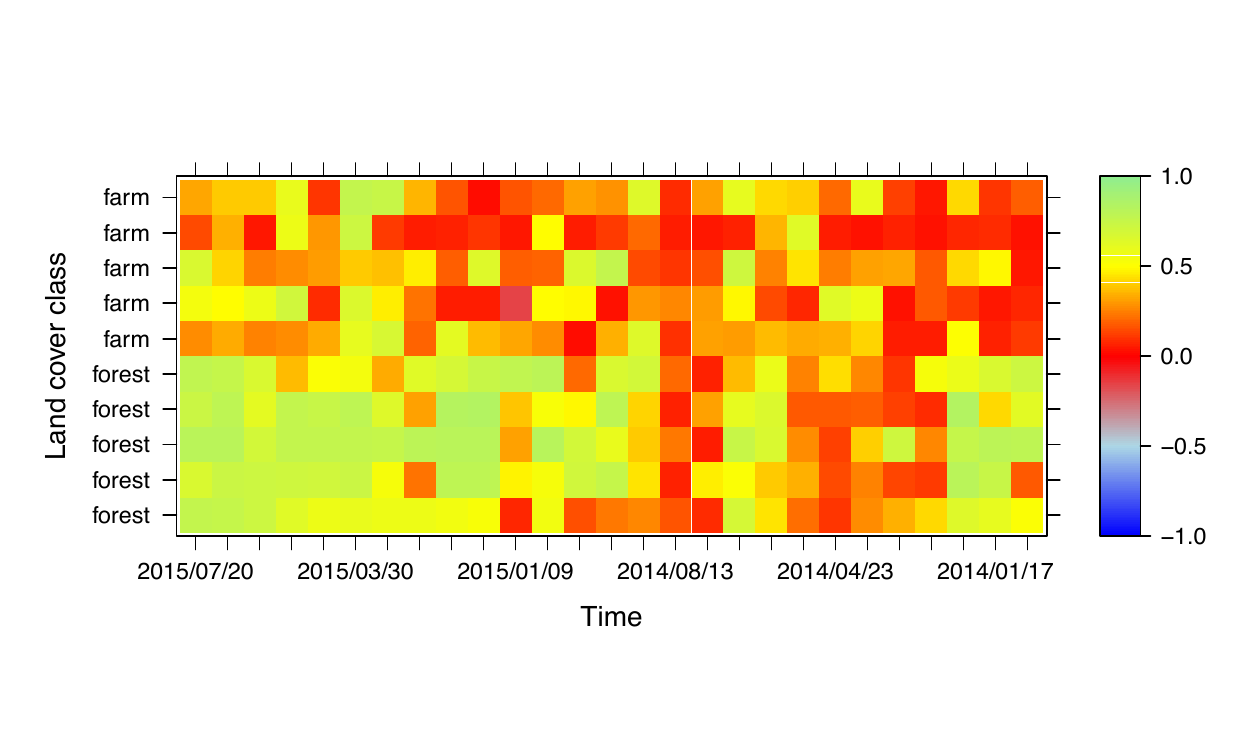}
		\caption{\small Illustration of the dataset.}
		\label{application_heatmap}
	\end{subfigure}
	\begin{subfigure}[b]{0.495\textwidth}
		\includegraphics[width=2.3in,height=2in]{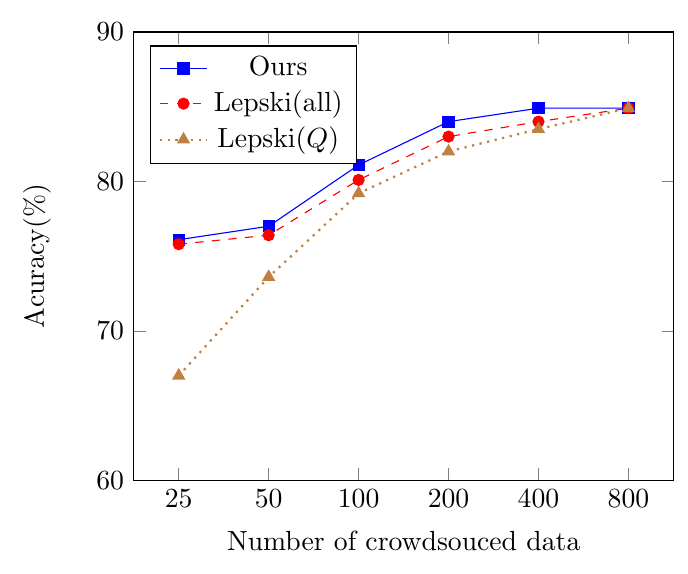}
		\caption{\small Accuracy of different classifiers.} 
		\label{application_results}
	\end{subfigure}
	\caption{\small (a) Illustration of the dataset. Each row represents one of a land cover class (farm or forest) and corresponding NDVI values of a pixel from remotely-sensed imagery in 2014-2015. (b) Accuracy of the three methods on the crowdsourced mapping data with different numbers of crowdsourced data involved. Blue: The proposed adaptive classifier. Red: Lepski's method using combined data. Brown: Lepski's method using only crowdsourced data.} \vspace{-15pt}
\end{figure}

Figure \ref{application_results} shows the accuracy of the three methods with different numbers of crowdsourced data involved. As more and more crowdsourced data  are used, the amount of information contained in the crowdsourced data  gradually increases, and the relative contribution from the accurately labeled data gradually decreases. The proposed adaptive classifier consistently outperforms the naive Lepski's method, especially when the number of the crowdsourced data is between 100 and 400, because in these cases the adaptive classifier can significantly increase the accuracy by better leveraging the information gained from both distributions.

\section{Multiple Source Distributions}
\label{Extension.sec}
	 
	We have so far focused on transfer learning with one source distribution $P$ and one target distribution $Q$. In practice, data may be generated from more than one source distribution. In this section, we generalize our methods to treat transfer learning in the setting  where multiple source distributions are available.
	
	We consider the following model where there are several source distributions with different relative signal exponents with respect to the target distribution $Q$. Suppose there are $n_{P_1},  ..., n_{P_m},$ and  $n_Q$ i.i.d data points generated from the source distributions $P_1, ..., P_m,$ and the target distribution $Q$ respectively.
	\baligns
	(X_1^{P_1}, Y_1^{P_1}),  ..., (X_{n_{P_1}}^{P_1}, Y_{n_{P_1}}^{P_1}) &\iid P_1  \\
	\vdots\\
	(X_1^{P_m}, Y_1^{P_m}), ..., (X_{n_{P_m}}^{P_m}, Y_{n_{P_m}}^{P_m}) &\iid P_m  \\
	(X_1^Q, Y_1^Q),  ..., (X_{n_Q}^Q, Y_{n_Q}^Q) &\iid Q
	\ealigns
and all the samples are independent. The goal is to make classification under the target distribution $Q$. Similar as before, it is intuitively clear that how useful the data from the source distributions $P_i$, $i\in[m]$, to the classification task under $Q$ depends on the relationship between each $P_i$ and $Q$. The definition of the relative signal exponent needs to be extended to accommodate the multiple source distributions. It is natural to consider the situation where each  source distribution $P_i$ and the target distribution $Q$ have a relative signal exponent. This motivates the following definition of the vectorized relative signal exponent when there are multiple source distributions.
 	 
	\begin{definition}
	Suppose the distributions $P_1, ..., P_m$, and $Q$ are supported on $\reals^d \times \{ 0,1 \}$. Define the class $\Gamma(\bm{\gamma}, C_{\bm{\gamma}})$ with  the relative signal exponent $\bm{\gamma} = (\gamma_1, ..., \gamma_m)  \in \reals_+^m$ and constants $C_{\bm{\gamma}} = (C_1, ..., C_m ) \in \reals_+^m$, is the set of distribution tuples $(P_1, ..., P_m, Q)$ that satisfy, for each $i \in [m]$, $(P_i, Q)$ belongs to the class $\Gamma(\gamma_i, C_i)$ with the relative signal exponent $\gamma_i$.
	\end{definition}
	
	Similar as in Section \ref{Model.sec}, adding the regularity conditions on $Q$ including the smoothness, margin assumption and strong density assumption, we define the parameter space in the multiple source distributions setting as follows:
		\baligns 
	\Pi(\alpha, & C_\alpha,  \beta, C_\beta, \bm{\gamma}, C_{\bm{\gamma}},  \mu, c_\mu, r_\mu) = 
	\{ (P_1, ..., P_m, Q) : (P_1, ..., P_m, Q) \in \Gamma(\bm{\gamma}, C_{\bm{\gamma}}),  \\ 
	&Q \in \mathcal{M}(\alpha, C_\alpha),  \eta_Q \in \mathcal{H}(\beta,C_\beta), (P_{i,X},Q_{X}) \in \mathcal{S}(\mu,c_\mu, r_\mu) \ \text{for all $i \in [m]$}  \} .
	\ealigns
	
	The above space will be simply denoted by $\Pi$ or $\Pi(\alpha, \beta, \bm{\gamma}, \mu)$ if there is no confusion.
	
	In this section we establish the minimax optimal rate of convergence and propose an adaptive classifier which simultaneously achieves the optimal rate of convergence within a logarithmic factor over a wide collection of the parameter spaces. The proofs are similar to those for Theorems \ref{ub}, \ref{lb} and \ref{adap_bound}  in the one source distribution setting. For reasons of space, we omit the proofs.
	
	\subsection{Minimax rate of convergence}
	
	We begin with the construction of a minimax rate optimal classifier $\hat f_{NN}$ in the case of multiple source distributions. The classifier is an extension of the two-sample weighted $K$-NN classifier given in Section  \ref{Minimax.sec}. It incorporates the information contained in the data drawn from the source distributions $P_i$, $i\in [m]$,  as well as the data drawn from the target distribution $Q$. The detailed steps are as follows.
	\begin{enumerate}
		\item[Step 1:] Compute the weights $w_{P_1}, ..., w_{P_m}$, and $w_Q$ by
		\baligns
		w_{P_i} &= (n_Q + \sum_{i=1}^m n_{P_i}^\frac{2\beta+d}{2\gamma_i\beta + d})^{-\frac{\gamma_i\beta}{2\beta+d}}, \quad \text{ for all } i\in [m]  \\
		w_Q &= (n_Q + \sum_{i=1}^m n_{P_i}^\frac{2\beta+d}{2\gamma_i\beta + d})^{-\frac{\beta}{2\beta+d}}.  
		\ealigns
		
		Compute the numbers of neighbors $k_{P_1},  ..., k_{P_m}, k_Q$ by
		\baligns 
		k_{P_i} &= \lfloor n_{P_i}(n_Q + \sum_{i=1}^m n_{P_i}^\frac{2\beta+d}{2\gamma_i\beta + d})^{-\frac{d}{2\beta+d}} \rfloor, \quad \text{ for all } i\in [m]  \\
		k_Q &= \lfloor n_Q(n_Q + \sum_{i=1}^m n_{P_i}^\frac{2\beta+d}{2\gamma_i\beta + d})^{-\frac{d}{2\beta+d}} \rfloor.  
		\ealigns
		
		\item[Step 2:] Define $X_{(j)}^{P_i}(x)$ to be the $j$-th nearest data point to $x$ among $X_1^{P_i},...,X_{n_{P_i}}^{P_i}$ and $Y^{P_i}_{(j)}(x)$ is its corresponding response (label). Likewise, let $X_{(j)}^Q(x)$ be the $j$-th data point to $x$ among $X_1^Q, ..., X_{n_Q}^Q$ and $Y^Q_{(j)}(x)$ is its corresponding response (label). Define the weighted $K$-NN estimator
		\[ \hat{\eta}_{NN}(x) = \frac{w_Q\sum_{j=1}^{k_Q}Y^Q_{(j)}(x) + \sum_{i=1}^m \p{w_{P_i}\sum_{j=1}^{k_{P_i}}Y^{P_i}_{(j)}(x)}}{w_Qk_Q + \sum_{i=1}^m w_{P_i}k_{P_i}} \]
where this estimator takes weighted average among $k_{P_i}$ nearest neighbors from those data points drawn from $P_i$, each with weight $w_{P_i}$, and $k_Q$ nearest neighbors from those data points drawn from $Q$, each with weight $w_Q$.
		
		\item[Step 3:] The final classifier is then defined as
		\[ \hf_{NN}(x) = \ind{\hat\eta_{NN}(x) > \frac{1}{2}}. \]
	\end{enumerate}

	 We now analyze the theoretical properties of the classifier $\hf_{NN}$.
	Theorem \ref{multi-ub} gives an upper bound for the excess risk $\mathcal{E}_Q(\hat f_{NN})$, while Theorem \ref{multi-lb} provides a matching lower bound on the excess risk for all estimators. These two theorems together establish the minimax rate of convergence and show that $\hf_{NN}$ attains the optimal rate. 
	
	\begin{theorem}[Upper Bound] \label{multi-ub}
		There exists a constant $C>0$ not depending on $n_P$ or $n_Q$,
		such that
		\[ \sup_{(P_1,...,P_m,Q) \in \Pi(\alpha,\beta,\bm{\gamma},\mu)} \E \mathcal{E}_Q(\hf_{NN}) \leq C(n_Q + \sum_{i=1}^m n_{P_i}^\frac{2\beta+d}{2\gamma_i\beta + d})^{-\frac{\beta(1+\alpha)}{2\beta+d}} . \]
	\end{theorem}

	\begin{theorem}[Lower Bound] 	\label{multi-lb}
		There exists a constant $c>0$ not depending on $n_P$ or $n_Q$,
		 such that 
		\[ \inf_{\hf} \sup_{(P_1,...,P_m,Q) \in \Pi(\alpha,\beta,\bm{\gamma},\mu)} \E \mathcal{E}_Q(\hf) \geq c(n_Q + \sum_{i=1}^m n_{P_i}^\frac{2\beta+d}{2\gamma_i\beta + d})^{-\frac{\beta(1+\alpha)}{2\beta+d}}. \]
	\end{theorem}

		Theorems \ref{multi-ub} and \ref{multi-lb} together yield the minimax optimal rate for transfer learning with multiple source distributions:
		\beq \inf_{\hat f}\sup_{(P_1,...,P_m,Q) \in \Pi(\alpha,\beta,\bm{\gamma},\mu)} \E \mathcal{E}_Q(\hf) \asymp (n_Q + \sum_{i=1}^m n_{P_i}^\frac{2\beta+d}{2\gamma_i\beta + d})^{-\frac{\beta(1+\alpha)}{2\beta+d}}. \label{multi-optimal-rate}\eeq
		
		As discussed in Section \ref{Minimax.sec}, here $n_{P_i}^\frac{2\beta+d}{2\gamma_i\beta + d}$ can be viewed as the effective sample size for data drawn from the source distribution $P_i$ when the information in this sample is transferred to help the classification task under the target distribution $Q$. Even when there are multiple source distributions, the transfer rate associated with $P_i$ remains to be $\frac{2\beta+d}{2\gamma_i\beta + d}$, which is not affected by the presence of the data drawn from the other source distributions.

	\subsection{Adaptive classifier}
	
	Again, in practice all the model parameters $\alpha,\beta,\bm{\gamma}$ and $\mu$ are typically unknown and the minimax classifier is not practical. It is desirable to construct a data driven classifier that does not rely on  the knowledge of the model parameters. A similar adaptive data-driven classifier can be developed. The detailed steps are summarized  in Algorithm \ref{algo_multi}.
	
	\begin{algorithm}[t]
		\caption{The Data Driven Classifier}
		\begin{algorithmic}
			\State \textbf{Input: } $x \in \supp(Q_X)$
			\For{$k = 1, ...,(n_Q + \sum_{i=1}^m n_{P_i} -1),(n_Q + \sum_{i=1}^m n_{P_i}) $}
			
			\State Find $k$ nearest neighbors $X_{(1)}(x), ..., X_{(k)}(x)$ to $x$ among all the covariates $\{X_j^Q : j \in [n_Q]\} \cup \bigcup_{i=1}^m\{ X_j^{P_i} : j \in [n_{P_i}] \}$. Suppose $k_{P_i}^{(k)}$ of them are from the distribution $P_i$, $i =1, \ldots, m$, and $k_Q^{(k)}$ of them are from $Q$. That is, the $k$ nearest neighbors are partitioned into $m+1$ parts according to which distribution they are drawn from.
			
			\State For each $i\in[m]$, Compute the nearest neighbor estimate for $\eta_{P_i}$ (If $k_{P_i}^{(k)}=0$, set $\hat\eta_{P_i}^{(k)} \gets \frac{1}{2}$ )
			\[\hat\eta_{P_i}^{(k)}(x) \gets \frac{1}{k_{P_i}^{(k)}}\sum_{j=1}^{k_{P_i}^{(k)}}Y^{P_i}_{(j)}(x)\] 
			and nearest neighbor estimate for $\eta_{Q}$ (If $k_Q^{(k)}=0$, set $\hat\eta_Q^{(k)} \gets \frac{1}{2}$ )
			\[ \hat\eta_Q^{(k)} \gets \frac{1}{k_Q^{(k)}}\sum_{i=1}^{k_Q^{(k)}} Y^P_{(i)}(x).\]
			
			Compute the positive signal-to-noise index
			\[ \hat r^{(k)}_+ \gets \ind{\eta_Q^{(k)} \geq \frac{1}{2}} k_Q^{(k)} \p{\eta_Q^{(k)}- \frac{1}{2}}^2 + \sum_{i=1}^m \ind{\eta_{P_i}^{(k)} \geq \frac{1}{2}}k_{P_i}^{(k)}\p{\eta_{P_i}^{(k)} - \frac{1}{2}}^2  \]
			and negative signal-to-noise index
			\[ \hat r^{(k)}_- \gets \ind{\eta_Q^{(k)} < \frac{1}{2}} k_Q^{(k)} \p{\eta_Q^{(k)}- \frac{1}{2}}^2 + \sum_{i=1}^m \ind{\eta_{P_i}^{(k)} < \frac{1}{2}}k_{P_i}^{(k)}\p{\eta_{P_i}^{(k)} - \frac{1}{2}}^2. \]
			
			\State Let $\hat r^{(k)}$ be the signal-to-noise ratio index calculated by
			\baligns \hat r^{(k)} \gets \max \left\{ \hat r^{(k)}_+, \hat r^{(k)}_-  \right\}.
			\ealigns
			\State Define the classifier 
			\[  \hat f^{(k)}(x) \gets \ind{\hat r^{(k)}_+ \geq \hat r^{(k)}_-}. \]
			\If {$\hat r^{(k)} > (d+3)\log(n_Q + \sum_{i=1}^m n_{P_i})$}
			\State	Stop and output $\hat f_a(x) \gets f^{(k)}(x)$.
			\EndIf
			\EndFor
			
			\State Output $\hat f_a(x) \gets \hf^{(k_m)}(x)$ where $k_m = \argmax_k \hat r^{(k)}$.
		\end{algorithmic}
		\label{algo_multi}
	\end{algorithm}

	It is clear from the construction that the classifier $\hf_a$ is a data-driven decision rule. Theorem \ref{multi_adap} below provides a theoretical guarantee for the excess risk of $\hf_a$ under the target distribution $Q$. In view of the optimal rate given in \eqref{multi-optimal-rate},  Theorem \ref{multi_adap} shows that  $\hf_a$ is adaptively nearly optimal over a wide range of parameter spaces.

	\begin{theorem} \label{multi_adap}
		Let $n=n_Q + \sum_{i=1}^m n_{P_i}$. There exists a constant $C>0$ such that for $\Pi= \Pi(\alpha,\beta,\bm{\gamma},\mu)$,		
		\[ \sup_{(P_1, ..., P_m,Q) \in \Pi} \E \mathcal{E}_Q(\hat{f}_a) \leq C \cdot  \p{\frac{n_Q}{\log n} + \sum_{i=1}^m\p{\frac{n_{P_i}}{\log n}}^{\frac{2\beta + d}{2\gamma_i\beta+d}} }^{-\frac{\beta(1+\alpha)}{2\beta+d}}. \]
	\end{theorem}
	
\clearpage


\section{Discussion}
\label{Discussion.sec}

We studied in this paper transfer learning under the posterior drift model and established the minimax rate of convergence. The optimal rate quantifies precisely the amount of information in the $P$-data that can be transferred to help classification under the target distribution $Q$.  A delicately designed data-driven adaptive classifier was also constructed and shown to be simultaneously within a log factor of the optimal rate over a large collection of parameter spaces. 

The results and techniques developed in this paper serve as a starting point for the theoretical analysis of other transfer learning problems. For example, in addition to classification, it is also of significant interest to characterize the relationship between the source distribution and the target distribution, so that the data from the source distribution $P$ can help in other statistical  problems under the target distribution $Q$. Examples include regression, hypothesis testing, and construction of confidence sets. We will investigate these transfer learning problems in the future.

Within the posterior drift framework of this paper, some of the technical assumptions given in the model formulation can be relaxed to a certain extent. For the smoothness parameter $\beta$, we focused on the case  $0<\beta\le 1$. It is possible to consider more general classes where $\beta$ can be larger than 1 with a carefully designed weighted $K$-NN method, as was introduced in \cite{samworth2012optimal}. Construction of such a weighted $K$-NN method is involved and we leave it as future work. For the assumptions on the support of the marginal distributions $P_X$ and $Q_X$, other than the strong density assumption, there are also weaker regularity conditions introduced in the literature. See, for example,  \cite{gadat2016, kpotufe2018marginal}. Similar results on  the minimax rate of convergence can be established under these different regularity conditions. The minimax and adaptive procedures should also be suitably modified.

\section{Proofs}
\label{Proof.sec}
	
	We prove Theorem \ref{ub} in this section and leave the proofs of other theorems and additional technical lemmas in the supplementary material \citep{CaiWei2019Supplement}. For readers' convenience, we begin by stating  Lepski's method for nonparametric classification in the conventional setting where there are only the $Q$-data.

	\subsection{Lepski's method}
	
	Algorithm \ref{lepski_algo} is one version of Lepski's method in nonparametric classification. We state the whole algorithm here for reference.
	
	\begin{algorithm}[H]
		\caption{Lepski's method \citep{lepski1997optimal}}
		\begin{algorithmic}
			\State \textbf{Input: } $n$ labeled samples $(X_i, Y_i) \in \reals^d \times \{0,1\}$, $i \in [n]$, and a point $x \in \reals^d$ to be classified.
			\State Set $\eta^-_0 \gets -\infty$ and $\eta^+_0 \gets +\infty$.
			\For{$k = 1, ...,(n_P+n_Q-1),(n_P+n_Q) $}
			\State Find $k$ nearest neighbor estimates $\hat\eta_k(x) = \frac{1}{k}\sum_{i=1}^k Y_{(i)}$, where $Y_{(i)}$ denote the label to $i$-th nearest covariates to $x$.
			\State Set $\eta^-_k \gets \eta^-_{k-1} \vee (\hat\eta_k(x) - \sqrt{\frac{d+3}{k}}\log n)$.
			\State Set $\eta^+_k \gets \eta^+_{k-1} \wedge (\hat\eta_k(x) + \sqrt{\frac{d+3}{k}}\log n)$.
			\If {$\eta^-_k > \frac{1}{2}$ or $\eta^+_k < \frac{1}{2}$}
				\State Stop and output $\hat f_L(x) \gets \ind{\hat\eta_k(x) \geq \frac{1}{2}}$.
			\EndIf
			\EndFor
			\State Output $\hat f_L(x) \gets \ind{\hat\eta_n(x) \geq \frac{1}{2}}$.
		\end{algorithmic}
		\label{lepski_algo}
	\end{algorithm}
	
	\subsection{Proof of Theorem \ref{ub}}
	
	First we define some new notations for convenience. In the proof, we use $\zeta_Q(x) = |\eta_Q(x)-\frac{1}{2}|$ and $\zeta_P(x) = |\eta_P(x)-\frac{1}{2}|$ to denote the signal strength. Let $\bar{Y}^Q_{(1:k_Q)}(x) := \frac{1}{k_Q}\sum_{i=1}^{k_Q} Y_{(i)}^Q(x)$ and $\bar{Y}^P_{(1:k_P)}(x) := \frac{1}{k_P}\sum_{i=1}^{k_P} Y_{(i)}^P(x)$ denote the average of $k_Q$ nearest neighbors in $Q-$data and $k_P$ nearest neighbors in $P-$data respectively. We will sometime omit $x$ in the notations such as $X_{(i)}^Q(x), X_{(i)}^P(x)$ if there is no confusion in the context. We also use the shorthand $X_{1:n_Q}^Q$ to denote the whole set of the $Q-$data covariates $\{X_1^Q, ..., X_{n_Q}^Q\}$. And similarly $X_{1:n_P}^P$ denotes $\{ X_1^P, ..., X_{n_P}^P \}$. We define $\E_{Y|X}(\cdot) = \E(\cdot|X_{1:n_Q}^Q, X_{1:n_P}^P)$ to denote the expectation conditional on the covariates of all data. Also, in following proofs we always assume $(P,Q) \in \Pi(\alpha,\beta,\gamma,\mu)$ so we will not state this assumption again in the lemmas.
	
	Before proving the theorem, we first state three useful lemmas. The first lemma $\ref{distb}$ provides a high probability uniform bound on the distance between any point and its $k-$th nearest neighbor.
	
	\begin{lemma}[$K$-NN Distance Bound] \label{distb}
		There exists a constant $C_D>0$ such that, with probability at least $1-C_D\frac{n_Q}{k_Q}\exp(-\frac{k_Q}{6})$, for all $x\in \Omega$,
	
		\beq \Vert X_{(k_Q)}^Q(x) - x \Vert \leq C_D(\frac{k_Q}{n_Q})^{\frac{1}{d}}. \label{distb_eq1}\eeq
		
		And with probability at least $1 - C_D\frac{n_P}{k_P}\exp(-\frac{k_P}{6})$, for all $x\in \Omega$,
	
		\beq \Vert X_{(k_P)}^P(x) - x \Vert \leq C_D(\frac{k_P}{n_P})^{\frac{1}{d}}. \label{distb_eq2} \eeq
		
	\end{lemma}
	
	Let $E_Q$ denote the event that Inequality \eqref{distb_eq1} holds for all $x \in \Omega$ and let $E_P$ denotes \eqref{distb_eq2} holds for all $x \in \Omega$. From lemma \ref{distb} we know
	\[ \P(E_Q) \geq 1-C_D\frac{n_Q}{k_Q}\exp(-\frac{k_Q}{6}) \]
	\[ \P(E_P) \geq 1 - C_D\frac{n_P}{k_P}\exp(-\frac{k_P}{6}). \]
	
	Lemma \ref{biasb} points out that when the signal is sufficiently strong, bias of $\bar Y^Q(x)$ and $\bar Y^P(x)$ will not be too large to overwhelm the signal.
	
	\begin{lemma}[Bias Bound] \label{biasb}
		There exist constants $c_b>0$ and $ C_b > 0$ such that:
		
		If a point $x\in\Omega$ satisfies $\zeta_Q(x) \geq 2C_\beta \Vert X^Q_{(k_Q)}(x) - x \Vert^\beta$, then we have
			\beq \E_{Y|X} (\bar{Y}^Q_{(1:k_Q)}(x)) - \frac{1}{2} \geq c_b\zeta_Q(x) \eqif f^\ast(x) = 1, \label{bisb1'} \eeq
			\beq \E_{Y|X} (\bar{Y}^Q_{(1:k_Q)}(x)) - \frac{1}{2} \leq -c_b\zeta_Q(x) \eqif f^\ast(x) = 0. \label{bisb2'} \eeq
			
		If a point $x\in\Omega$ satisfies $\zeta_Q(x) \geq 2C_\beta \Vert X^P_{(k_P)}(x) - x \Vert^\beta$, then we have 
			\beq \E_{Y|X} (\bar{Y}^P_{(1:k_P)}(x)) - \frac{1}{2} \geq c_b\zeta_Q(x)^\gamma \eqif f^\ast(x) = 1, \label{bisb3'} \eeq
			\beq \E_{Y|X} (\bar{Y}^P_{(1:k_P)}(x)) - \frac{1}{2} \leq -c_b\zeta_Q(x)^\gamma \eqif f^\ast(x) = 0. \label{bisb4'} \eeq

		As a consequence, if a point $x\in\Omega$ satisfies $\zeta_Q(x) \geq C_b(\max\{\frac{k_Q}{n_Q},\frac{k_P}{n_P}\})^{\frac{\beta}{d}}$, then 
		\begin{itemize}
			\item Under the event $E_Q$, we have 
			\beq \E_{Y|X} (\bar{Y}^Q_{(1:k_Q)}(x)) - \frac{1}{2} \geq c_b\zeta_Q(x) \eqif f^\ast(x) = 1, \label{bisb1} \eeq
			\beq \E_{Y|X} (\bar{Y}^Q_{(1:k_Q)}(x)) - \frac{1}{2} \leq -c_b\zeta_Q(x) \eqif f^\ast(x) = 0. \label{bisb2} \eeq
			
			\item Under the envet $E_P$, we have
			\beq \E_{Y|X} (\bar{Y}^P_{(1:k_P)}(x)) - \frac{1}{2} \geq c_b\zeta_Q(x)^\gamma \eqif f^\ast(x) = 1, \label{bisb3} \eeq
			\beq \E_{Y|X} (\bar{Y}^P_{(1:k_P)}(x)) - \frac{1}{2} \leq -c_b\zeta_Q(x)^\gamma \eqif f^\ast(x) = 0. \label{bisb4} \eeq
		\end{itemize}

	\end{lemma}
	
	Lemma \ref{tailb} gives a bound on the probability of misclassification at certain covariates $x$. 
	\begin{lemma}[Misclassification Bound] \label{tailb}
		Let $C_b$ and $c_b$ be the constants defined in Lemma \ref{biasb}. If $\zeta_Q(x) \geq C_b(\max\{\frac{k_Q}{n_Q},\frac{k_P}{n_P}\})^{\frac{\beta}{d}}$, then
		\begin{itemize}
			\item Under the event $E_Q$, we have
			\[ \P_{Y|X}(\hf_{NN}(x) \neq f^\ast_Q(x)) \leq \exp\p{ -2\frac{[(c_bw_Qk_Q\zeta_Q(x) - w_Pk_P)\vee 0]^2}{k_Pw_P^2 + k_Qw_Q^2} }. \]
			
			\item Under the event $E_P$, we have
			\[ \P_{Y|X}(\hf_{NN}(x) \neq f^\ast_Q(x)) \leq \exp\p{ -2\frac{[(c_bw_Pk_P\zeta_Q(x)^\gamma - w_Qk_Q)\vee 0]^2}{k_Pw_P^2 + k_Qw_Q^2} }. \]
			
			\item Under the event $E_P\cap \E_Q$, we have
			\[ \P_{Y|X}(\hf_{NN}(x) \neq f^\ast_Q(x)) \leq \exp\p{ -2c_b^2\frac{(w_Pk_P\zeta_Q(x)^\gamma + w_Qk_Q\zeta_Q(x))^2}{k_Pw_P^2 + k_Qw_Q^2} }. \]
		\end{itemize}

	\end{lemma}

	Given the three lemmas above, the remain proof generally follows the proof of Lemma 3.1 in \cite{audibert2007fast}. Let $\delta = (n_P^{\frac{2\beta+d}{2\gamma\beta+d}}+n_Q)^{-\frac{\beta}{2\beta+d}}$.  When $w_P,w_Q,k_P,k_Q$ are equal to the values defined in Theorem \ref{ub}, we have
	\beq w_Q = \delta, w_P = \delta^\gamma, k_Q = \lfloor n_Q\delta^{\frac{d}{\beta}} \rfloor, k_P = \lfloor n_P\delta^{\frac{d}{\beta}} \rfloor.  \label{def_paras} \eeq
	
	We will approximate $k_Q = n_Q\delta^{\frac{d}{\beta}}$ and $k_P = n_P\delta^{\frac{d}{\beta}}$ in the following proof because one can easily prove such an approximation will only result in changing the constant factor of the upper bound.
	
	Now we state another lemma which shows a local misclassification bound with high probability when the parameters in the weighted K-NN estimator are properly chosen.
	
	\begin{lemma} \label{misclassb}
		Using $w_P,w_Q,k_P,k_Q$ defined in theorem \ref{ub} to construct a weighted K-NN estimator $\hat f_{NN}$. Then there exist constants $c_1, C_1>0$ such that, with probability at least $1- 2(n_P^{\frac{2\beta + d}{2\gamma\beta+d}} + n_Q)^{-\frac{\beta(1+\alpha)}{2\beta+d}}$, for all $x$ we have
		\beq \P_{Y|X}(\hf_{NN}(x) \neq f^\ast_Q(x)) \leq C_1\exp \p{-c_1(\frac{\zeta_Q(x)}{\delta})^{1\wedge\gamma}}. \label{c_1} \eeq
	\end{lemma}
	
	Denote  by $E_0$ the event that inequality $\eqref{c_1}$ holds for all $x$. Lemma \ref{misclassb} implies
	\[ \P(E_0) \geq 1- 2(n_P^{\frac{2\beta + d}{2\gamma\beta+d}} + n_Q)^{-\frac{\beta(1+\alpha)}{2\beta+d}}.  \]
	
	Consider the disjoint sets $\A_j \subset \Omega, j=0,1,2,...$ defined as
	\[ \A_0 := \{ x \in \Omega: 0<\zeta_Q(x) \leq \delta \},  \]
	\[ \A_j := \{ x \in \Omega: 2^{j-1}\delta<\zeta_Q(x) \leq 2^j\delta \} \text{ for } j\geq1.  \]
	
	Note that by the margin assumption, for all $j$,
	\[ Q_X(A_j) \leq Q_X(|\eta_Q - \frac{1}{2}| \leq 2^j\delta) \leq C_\alpha 2^{\alpha j}\delta^\alpha. \]
	
	Based on the partition $A_0, A_1, ...$ and the dual representation of $\mathcal{E}_Q(\hf)$ shown in \eqref{dual_rep}, we have a decomposition of $\mathcal{E}_Q(\hf_{NN})$:
	\baligns
		\mathcal{E}_Q(\hf_{NN}) &= 2\E_{X\sim Q_X}(|\eta_Q(X) - \frac{1}{2}|\ind{ \hf_{NN}(X) \neq f_Q^\ast(X)}) \\
		&= 2\sum_{j=0}^\infty \E_{X\sim Q_X}(\zeta_Q(X)\ind{ \hf_{NN}(X) \neq f_Q^\ast(X)}\ind{X\in A_j}).
	\ealigns
	
	For $j=0$ we have
	\[  \E_{X\sim Q_X}(\zeta_Q(X)\ind{ \hf_{NN}(X) \neq f_Q^\ast(X)}\ind{X\in A_0}) \leq \delta \cdot Q_X(A_0) \leq C_\alpha \delta^{\alpha+1}.  \]
	
	Under the event $E_0$, $2^{j-1}\delta < \zeta(x) \leq 2^j\delta$ for $x \in A_j$ and $j>1$. Inequality \eqref{c_1} now yields
	\baligns
		\E_{Y|X}\E_{X\sim Q_X}&(\zeta_Q(X)\ind{ \hf_{NN}(X) \neq f_Q^\ast(X)}\ind{X\in A_j}) \\ &= \E_{X\sim Q_X}(\zeta_Q(X) \P_{Y|X}(\hf_{NN}(X) \neq f_Q^\ast(X))\ind{X\in A_j})\\
		&\leq 2^j\delta \cdot C_1\exp(-c_1\cdot2^{(j-1)\cdot (1\wedge \gamma)})\cdot Q_X(A_j) \\
		&\leq C_\alpha C_1[2^{(1+\alpha)j}\exp(-c_1\cdot2^{(j-1)\cdot (1\wedge \gamma)})]\delta^{\alpha+1}.
	\ealigns
		
	Combine those summands together we can obtain 
	\baligns \E_{Y|X}\mathcal{E}_Q(\hf_{NN}) &= 2\sum_{j=0}^\infty \E_{Y|X}\E_{X\sim Q_X}(\zeta_Q(X)\ind{ \hf_{NN}(X) \neq f_Q^\ast(X)}\ind{X\in A_j}) \\
	&\leq 2C_\alpha\p{1+C_1\sum_{j=0}^\infty [2^{(1+\alpha)j}\exp(-c_1\cdot2^{(k-1)\cdot (1\wedge \gamma)})]}\delta^{1+\alpha}  \\
	&\leq C\delta^{1+\alpha}.
	\ealigns
	
	Because summation $\sum_{j=0}^\infty [2^{(1+\alpha)j}\exp(-c_1\cdot2^{(k-1)\cdot (1\wedge \gamma)})]$ converges when $\gamma > 0$.
	
	Finally, it follows from Lemma \ref{misclassb} that
	\[ \P(E_0^c) \leq 2(n_P^{\frac{2\beta + d}{2\gamma\beta+d}} + n_Q)^{-\frac{\beta(1+\alpha)}{2\beta+d}}.  \]
	
	Applying the trivial bound $\mathcal{E}_Q(\hf_{NN}) \leq 1$ when $E_0^c$ occurs, we have
	\baligns \E \mathcal{E}_Q(\hf_{NN}) &= \E(\E_{Y|X}\mathcal{E}_Q(\hf_{NN})) \\
	&\leq \E(\E_{Y|X}\mathcal{E}_Q(\hf_{NN}) | E_0)\P(E_0) + \E(\E_{Y|X}\mathcal{E}_Q(\hf_{NN}) | E_0^c)\P(E_0^c) \\ 
	&\leq C(n_P^{\frac{2\beta + d}{2\gamma\beta+d}} + n_Q)^{-\frac{\beta(1+\alpha)}{2\beta+d}} \cdot 1 + 1\cdot 2(n_P^{\frac{2\beta + d}{2\gamma\beta+d}} + n_Q)^{-\frac{\beta(1+\alpha)}{2\beta+d}} \\
	&= (C+2)(n_P^{\frac{2\beta + d}{2\gamma\beta+d}} + n_Q)^{-\frac{\beta(1+\alpha)}{2\beta+d}}.
	\ealigns
	
	\qed 

\begin{supplement}[id=supplemnt]
	\sname{Supplement A}
	\stitle{Supplement to ``Transfer Learning for Nonparametric Classification: Minimax Rate and Adaptive Classifier"}
	\slink[doi]{url to be specified}
	\sdescription{In this supplementary material, we provide proofs for Theorems \ref{lb} and \ref{adap_bound}, and proofs for technical lemmas \ref{distb}, \ref{biasb}, \ref{tailb}, and \ref{misclassb}. The proofs of Theorems \ref{multi-ub}, \ref{multi-lb}, and \ref{multi_adap} are similar and thus omitted.}
\end{supplement}

\bibliographystyle{apalike}
\bibliography{Transfer-Learning-Ref}

\newpage	
%


\begin{frontmatter}
\title{Supplement to ``Transfer Learning for Nonparametric Classification: Minimax Rate and Adaptive Classifier"\thanksref{T1}}
\runtitle{Supplement to ``Transfer Learning"}
\thankstext{T1}{The research was supported in part by NSF Grant DMS-1712735 and NIH grants R01-GM129781 and R01-GM123056.}
\begin{aug}
\author{\fnms{T. Tony} \snm{Cai}\ead[label=e1]{tcai@wharton.upenn.edu}},
\and
\author{\fnms{Hongji} \snm{Wei}\ead[label=e2]{hongjiw@wharton.upenn.edu}
\ead[label=u1,url]{URL: http://www-stat.wharton.upenn.edu/$\sim$tcai/}}
\runauthor{T. T. Cai and H. Wei}
\affiliation{University of Pennsylvania}
\address{DEPARTMENT OF STATISTICS\\
THE WHARTON SCHOOL\\
UNIVERSITY OF PENNSYLVANIA\\
PHILADELPHIA, PENNSYLVANIA 19104\\
USA\\
\printead{e1}\\
\phantom{E-mail:\ }\printead*{e2}\\
\printead*{u1}\phantom{URL:\ }\\
}

\end{aug}

\begin{abstract}
We present in this supplement the detailed proofs of Theorems 2 and 3 in the paper ``Transfer Learning for Nonparametric Classification: Minimax Rate and Adaptive Classifier". We also prove the technical lemmas 1, 2, 3, and 4.
\end{abstract}

\end{frontmatter}

\setcounter{page}{1}
\setcounter{section}{0}

\section{Proof of Auxiliary lemmas of Theorem 5}
	\subsection{Proof of Lemma \ref{distb}}
	We only prove Equation \eqref{distb_eq1} as Equation \eqref{distb_eq2} can be proved in a similar way.
	
	Let $B(x,r) = \{ x': \Vert x'-x \Vert \leq r \}$ denote a ball centered at $x$ with radius $r$. Recall by our assumptions on $(P,Q)$, the marginal desity $\frac{dQ_X}{d\lambda}(x)$ is lower bounded by $\mu_-$ when $x \in \Omega$. Therefore for any $x \in \Omega$, $r<r_\mu$,
	\beq Q(X \in B(x,r)) = \int_{B(x,r)\cap\Omega} \frac{dQ_X}{d\lambda}(y) dy \geq \mu_-\lambda(B(x,r)\cap\Omega) \geq c_\mu \mu_-\pi_dr^d  \label{biasb1} \eeq
	where $\pi_d = \lambda(B(0,1))$ is the volume of $d$ dimension unit ball.
	
	When $\frac{k_Q}{n_Q} \geq \min\{\frac{1}{4},\frac{c_\mu\mu_-\pi_dr_\mu^d}{2}\}$, we can set $C_D$ large enough to make equality $\eqref{distb_eq1}$ match the trivial bound $\Vert X_{(k_Q)}^Q(x) - x \Vert \leq \sqrt{d}$. Therefore in the following proof we focus on the case when $\frac{k_Q}{n_Q} \leq \min\{\frac{1}{4},\frac{c_\mu\mu_-\pi_dr_\mu^d}{2}\}$.
	
	Set $r_0 = (\frac{2k_Q}{c_\mu\mu_-\pi_dn_Q})^{\frac{1}{d}}$, by $\frac{k_Q}{n_Q} < \frac{c_\mu\mu_-\pi_dr_\mu^d}{2}$ we have $r_0 < r_\mu$, thus from $\eqref{biasb1}$ we have for any $x\in \Omega$,
	\[  Q(X \in B(x,r_0)) \geq \frac{2k_Q}{n_Q}. \]
	
	Note that $\ind{X_i^Q \in B(x,r_0)}$ are i.i.d Bernuolli variables with mean $Q(X \in B(x,r_0)) \geq \frac{2k_Q}{n_Q}$. Let $S_n(x) = \sum_{i=1}^{n_Q} \ind{X_i^Q \in B(x,r_0)}$ represent the number of $Q-$data whose covariates fall into $B(x,r_0)$. Set $W \sim Binomial(n_Q, \frac{2k_Q}{n_Q})$, by Berstein inequality,
	\beq \P(S_n(x) < k_Q) \leq \P(W < k_Q ) = \P(W-\E W < -k_Q) \leq \exp(-\frac{k_Q^2}{2(2k_Q + k_Q)}) = \exp(-\frac{k_Q}{6}). \label{dist_proof1} \eeq 
	
	Inequality \eqref{dist_proof1} implies that the probability of $S_n(x) < k_Q$ is small for any given $x \in \Omega$. To get a union bound, we need to apply this bound on a set of balls covering the whole support. 
	
	Let $\mathcal{B} \subset \Omega$ be a finite set such that
	\[ \Omega \subset \bigcup_{x\in\mathcal{B}} B(x,r_0).  \]
	
	Note that $\Omega \subset [0,1]^d$, thus by easy construction we can find a feasible set $\mathcal{B}$ with $|\mathcal{B}| \leq Cr_0^{-d}$. This leads to a union bound
	\balign 
		\P(\exists x\in \mathcal{B}, S_n(x) < k_Q) \leq \sum_{x\in \mathcal{B}} \P(S_n(x) < k_Q) \leq Cr_0^{-d}\exp(-\frac{k_Q}{6}).
	\ealign
	
	For any $x \in \Omega$, there exist $x' \in \mathcal{B}$ such that $x \in B(x',r_0)$. Under the event $E_2 := \{ \forall x \in \mathcal{B}, S_n(x) \geq k_Q \}$, there are at least $k_Q$ $Q-$data covariates among $X_1^Q, ..., X_n^Q$ in the ball $B(x',r_0)$. All these points have distance to $x$ no larger than $2r_0$. Thus we have $\Vert X^Q_{(k_Q)}(x) - x \Vert \leq 2r_0$. Therefore,
	\[ \P(\forall x, \Vert X^Q_{(k_Q)}(x) - x \Vert \leq 2r_0) \geq \P(E_2) \geq 1-Cr_0^d\exp(-\frac{k_Q}{6}).  \]
	
	Plug in $r_0 = (\frac{2k_Q}{c_\mu\mu_-\pi_dn_Q})^{\frac{1}{d}}$ we can conclude that with probability at least $1 - \frac{2C}{c_\mu\mu_-\pi_d}\frac{n_Q}{k_Q}\exp(-\frac{k_Q}{6})$, we have $\Vert X^Q_{(k_Q)}(x) - x \Vert \leq 2(\frac{2}{c_\mu\mu_-\pi_d})^\frac{1}{d}(\frac{k_Q}{n_Q})^\frac{1}{d}$. Set $C_D = \max(\frac{2C}{c_\mu\mu_-\pi_d}, 2(\frac{2}{c_\mu\mu_-\pi_d})^\frac{1}{d})$ we can obtain the desired bound.
	
	Apply the similar result on $P-$data the lemma can be proved. \qed
	
	\subsection{Proof of Lemma \ref{biasb}}
	
	\textbf{1. Proof of \eqref{bisb1'}, \eqref{bisb2'}}
	
	Note that $\E_{Y|X}(Y^Q_{(i)}(x)) = \eta_Q(X_{(i)}^Q(x))$. When $\zeta_Q(x) \geq 2C_\beta \Vert X^Q_{(k_Q)}(x) - x \Vert^\beta$ we have
	\baligns
		\left|\E_{Y|X}\p{\bar Y^Q_{(1:k_Q)}(x) } - \eta_Q(x)\right|
		&\leq \frac{1}{k_Q}\sum_{i=1}^{k_Q}\left|\E_{Y|X}\p{Y_{(i)}^Q(x)} - \eta_Q(x)\right| \\
		&= \frac{1}{k_Q}\sum_{i=1}^{k_Q} \left| \eta_Q(X_{(i)}^Q(x)) - \eta_Q(x)\right| \\ 
		&\leq \frac{1}{k_Q}\sum_{i=1}^{k_Q}  C_\beta\Vert X_{(i)}^Q(x) - x \Vert^\beta \\
		&\leq C_\beta \Vert X_{(k_Q)}^Q(x) - x \Vert^\beta \\
		&\leq \frac{1}{2}\zeta_Q(x).
	\ealigns
	
	When $f^\ast = 1$ we have $\eta_Q(x) - \frac{1}{2} = \zeta_Q(x)$, thus
	\[ \E_{Y|X}\p{\bar Y^Q_{(1:k_Q)}(x) } - \frac{1}{2} \geq (\eta_Q(x) - \frac{1}{2}) - 	\left|\E_{Y|X}\p{\bar Y^Q_{(1:k_Q)}(x) } - \eta_Q(x)\right| \geq \frac{1}{2}\zeta_Q(x). \]
	
	When $f^\ast = 0$ we have $\eta_Q(x) - \frac{1}{2} = -\zeta_Q(x)$, thus
	\[ \E_{Y|X}\p{\bar Y^Q_{(1:k_Q)}(x) } - \frac{1}{2} \leq (\eta_Q(x) - \frac{1}{2}) + \left|\E_{Y|X}\p{\bar Y^Q_{(1:k_Q)}(x) } - \eta_Q(x)\right| \leq -\frac{1}{2}\zeta_Q(x). \]
	
	Therefore \eqref{bisb1'}, \eqref{bisb2'} holds as long as $c_b < \frac{1}{2}$.
	
	\textbf{2. Proof of \eqref{bisb3'}, \eqref{bisb4'}}
	
	Note that $\E(Y^P_{(i)}(x)|X_{1:n_P}^P) = \eta_P(X_{(i)}^P(x))$. When $f^\ast(x) = 1$ and $\zeta_Q(x) \geq 2C_\beta \Vert X^P_{(k_P)}(x) - x \Vert^\beta$ we have
	\baligns
	\E_{Y|X}\p{\bar Y^P_{(1:k_P)}(x)}  - \frac{1}{2}
	&\geq \frac{1}{k_P}\sum_{i=1}^{k_Q} \E_{Y|X}\p{Y_{(i)}^P(x) - \frac{1}{2}} \\
	&= \frac{1}{k_P}\sum_{i=1}^{k_P} \p{ \eta_P(X_{(i)}^P(x)) - \frac{1}{2} } \\ 
	&\geq \frac{1}{k_P}\sum_{i=1}^{k_P} C_\gamma\p{ \eta_Q(X_{(i)}^P(x)) - \frac{1}{2} }^\gamma \\ 
	&\geq \frac{1}{k_P}\sum_{i=1}^{k_P} C_\gamma\p{ \max\{ \eta_Q(x) - \frac{1}{2} - C_\beta \Vert (X_{(i)}^P(x) - x \Vert^\beta , 0\} }^\gamma \\ 
	&\geq  C_\gamma\p{ \max\{ \zeta_Q(x) - C_\beta \Vert (X_{(k_P)}^P(x) - x \Vert^\beta , 0\} }^\gamma \\
	&\geq  C_\gamma\p{\frac{1}{2}\zeta_Q(x) }^\gamma.
	\ealigns
	
	Similarly when $f^\ast(x) = 0$ and $\zeta_Q(x) \geq 2C_\beta \Vert X^P_{(k_P)}(x) - x \Vert$ we can obtain 
	\[ \E_{Y|X}\p{\bar Y^P(x)}  - \frac{1}{2} \leq -C_\gamma\p{\frac{1}{2}\zeta_Q(x) }^\gamma.  \]
	
	Therefore \eqref{bisb3'}, \eqref{bisb4'} holds as long as $c_b < \frac{C_\gamma}{2^\gamma}$.
	
	\textbf{3. Proof of \eqref{bisb1}, \eqref{bisb2}, \eqref{bisb3}, \eqref{bisb4}}
	
	Here we set $C_b = 2C_D^\beta C_\beta$.
	
	Under the event $\E_Q$, from \eqref{distb_eq1} we know $\Vert X_{(k_Q)}^Q(x) - x \Vert \leq C_D(\frac{k_Q}{n_Q})^{\frac{1}{d}}$. So $\zeta_Q(x) \geq C_b(\max\{\frac{k_Q}{n_Q},\frac{k_P}{n_P}\})^{\frac{\beta}{d}}$ implies
	
	\[ \zeta_Q(x) \geq 2C_\beta C_D^\beta (\frac{k_Q}{n_Q})^\frac{\beta}{d} \geq 2C_\beta \Vert X^Q_{(k_Q)} - x \Vert^\beta. \]
	
	The above inequality is exactly the condition such that \eqref{bisb1'} and \eqref{bisb2'} holds. Therefore we have \eqref{bisb1} and \eqref{bisb2}.
	
	Under the event $\E_P$, from \eqref{distb_eq2} we know $\Vert X_{(k_P)}^P(x) - x \Vert \leq C_D(\frac{k_P}{n_P})^{\frac{1}{d}}$. So $\zeta_Q(x) \geq C_b(\max\{\frac{k_Q}{n_Q},\frac{k_P}{n_P}\})^{\frac{\beta}{d}}$ implies
	 
	\[ \zeta_Q(x) \geq 2C_\beta C_D^\beta (\frac{k_P}{n_P})^\frac{\beta}{d} \geq 2C_\beta \Vert X^P_{(k_P)} - x \Vert^\beta. \]
	
	The above inequality is exactly the condition such that \eqref{bisb3'} and \eqref{bisb4'} holds. Therefore we have \eqref{bisb3} and \eqref{bisb4}.
	
	\qed
	
	\subsection{Proof of Lemma \ref{tailb}}

	We will treat $x$ as a fixed point during the following proof. Without loss of generality, we assume $f^\ast(x) = 1$. The case $f^\ast(x) = 0$ can be proved by a similar way. 
	
	The following proofs are derived under the condition $\zeta_Q(x) \geq C_b(\max\{\frac{k_Q}{n_Q},\frac{k_P}{n_P}\})^{\frac{\beta}{d}}$. From lemma $\ref{biasb}$ we know under $E_Q$ \eqref{bisb1} and \eqref{bisb2} hold, under $E_P$ \eqref{bisb3} and \eqref{bisb4} hold.
	
	Note that 
	\baligns 
		\E_{Y|X}(\hat\eta_{NN}(x)- \frac{1}{2})\vee 0 &= \E_{Y|X}(\frac{w_Pk_P\bar Y^P_{(1:k_P)} + w_Qk_Q\bar Y^Q_{(1:k_Q)}}{w_Pk_P + w_Qk_Q}- \frac{1}{2})\vee 0 \\ 
		&= \frac{w_Pk_P\p{\E_{Y|X} (\bar Y^P_{(1:k_P)}) - \frac{1}{2}} + w_Qk_Q\p{ \E_{Y|X}(\bar Y^Q_{(1:k_Q)}) - \frac{1}{2}}}{w_Pk_P + w_Qk_Q} \vee 0. 
	\ealigns
	
	So under $E_Q$, from \eqref{bisb1} and $\E_{Y|X}(\bar Y^P_{(1:k_P)}) - \frac{1}{2} > -1$ we have
	\beq
		\E_{Y|X}(\hat\eta_{NN}(x)- \frac{1}{2})\vee 0 \geq \frac{c_bw_Qk_Q\zeta_Q(x) - w_Pk_P}{w_Pk_P + w_Qk_Q} \vee 0.
		\label{taib1}
	\eeq
	
	Under $E_P$, from \eqref{bisb3} and $\E_{Y|X}(\bar Y^Q_{(1:k_Q)}) - \frac{1}{2} > -1$ we have
	\beq
	\E_{Y|X}(\hat\eta_{NN}(x)- \frac{1}{2})\vee 0 \geq \frac{c_bw_Pk_P\zeta_Q(x)^\gamma - w_Qk_Q}{w_Pk_P + w_Qk_Q} \vee 0.
	\label{taib2}
	\eeq
	
	Under $E_P\cap E_Q$, from \eqref{bisb1}, \eqref{bisb3} we have
	\beq
	\E_{Y|X}(\hat\eta_{NN}(x)- \frac{1}{2})\vee 0 \geq c_b\frac{w_Pk_P\zeta_Q(x)^\gamma + w_Qk_Q\zeta_Q(x)}{w_Pk_P + w_Qk_Q}.
	\label{taib3}
	\eeq
	
	Moreover, observe the formula
	\[ \hat\eta_{NN}(x) - \frac{1}{2} = \sum_{i=1}^{k_P} \frac{w_P}{w_Pk_P + w_Qk_Q}Y_{(i)}^P(x) + \sum_{i=1}^{k_Q} \frac{w_Q}{w_Pk_P + w_Qk_Q}Y_{(i)}^Q(x) - \frac{1}{2}.  \]
	
	Note that $Y_{(i)}^P(x)\in\{0,1\}$ for all $i\in [n_P]$, changing any single entry $Y_{(i)}^P(x)$ will result in changing $\hat\eta_{NN}(x) - \frac{1}{2}$ at most $\frac{w_P}{w_Pk_P+w_Qk_Q}$. For the same reason, changing any single entry $Y_{(i)}^Q(x)$ will result in changing $\hat\eta_{NN}(x) - \frac{1}{2}$ at most $\frac{w_Q}{w_Pk_P+w_Qk_Q}$. Condition on $X_{1:n_P}^P \cup X_{1:n_Q}^Q$, $Y_{(1)}^P, ..., Y_{(n_P)}^P$ and $Y_{(1)}^Q, ..., Y_{(n_Q)}^Q$ are all independent, thus by McDiarmid's inequality,
	
	\baligns
		\P_{Y|X}(\hf_{NN}(x) \neq f^\ast_Q(x)) &= \P_{Y|X}(\hat\eta_{NN}(x) - \frac{1}{2}\leq0 ) \\
		&= \P_{Y|X}\p{(\hat\eta_{NN}(x) - \frac{1}{2}) - \E_{Y|X}(\hat\eta_{NN}(x) - \frac{1}{2}) \leq -\E_{Y|X}(\hat\eta_{NN}(x) - \frac{1}{2})} \\
		&\leq \exp\p{ -\frac{  2(\E_{Y|X}(\hat\eta_{NN}(x) - \frac{1}{2}) \vee 0)^2  } {  k_P(\frac{w_P}{w_Pk_P+w_Qk_Q})^2 + k_Q(\frac{w_Q}{w_Pk_P+w_Qk_Q})^2  } }.
	\ealigns 
	
	Plug in \eqref{taib1}, \eqref{taib2} or \eqref{taib3} we can obtain the desired bounds stated in the lemma under event $E_Q$, $E_P$ or $E_P\cap E_Q$ respectively.
	\qed
	
	\subsection{Proof of Lemma \ref{misclassb}}
	When $\zeta_Q(x) < C_b\delta$ the probability bound \eqref{c_1} is trivial.
	
	So from now on we assume $\zeta_Q(x) \geq C_b\delta$. Then with specific choices of $k_P, k_Q$ stated in \eqref{def_paras} we have
	\beq \zeta_Q(x) \geq C_b\delta \geq C_b(\max\{\frac{k_Q}{n_Q},\frac{k_P}{n_P}\})^{\frac{\beta}{d}}. \label{mcb1} \eeq
	
	For simplicity, denote $t = \frac{\zeta_Q(x)}{C_b\delta}$. We are going to discuss 4 cases depending on the values of $n_P$ and $n_Q$. In each case we will show \eqref{c_1} holds with desired high probability.
	
	\textbf{Case 1:} $\frac{n_Q}{k_Q}\exp(-\frac{k_Q}{6}) < \delta^{1+\alpha}$ and $\frac{n_P}{k_P}\exp(-\frac{k_P}{6}) < \delta^{1+\alpha}$
	
	Then based on \eqref{mcb1} in lemma $\ref{tailb}$, we have under event $E_P\cap E_Q$, given $C_b > 1$ large enough, we have for all $x\in\Omega$ that satisfies $\zeta_Q(x) \geq C_b\delta$,
	\baligns 
	\P_{Y|X}(\hf_{NN}(x) \neq f^\ast_Q(x)) &\leq \exp\p{ -2c_b^2\frac{(w_Pk_P\zeta_Q(x)^\gamma + w_Qk_Q\zeta_Q(x))^2}{k_Pw_P^2 + k_Qw_Q^2} } \\ 
	&= \exp \p{ -2c_b^2\frac{(\delta^\gamma \cdot n_P\delta^{\frac{d}{\beta}} \cdot (C_bt\delta)^\gamma + \delta \cdot n_Q\delta^{\frac{d}{\beta}} \cdot C_bt\delta)^2}{n_P\delta^{\frac{d}{\beta}}\cdot \delta^{2\gamma} + n_Q\delta^{\frac{d}{\beta}} \cdot \delta^2   }  }	\\
	& = \exp \p{  -2c_b^2\frac{[(C_bt)^\gamma n_P\delta^{2\gamma+\frac{d}{\beta}} + (C_bt)n_Q\delta^{2+\frac{d}{\beta}}]^2}{n_P\delta^{2\gamma+\frac{d}{\beta}} + n_Q\delta^{2+\frac{d}{\beta}}   }  } \\
	& \leq \exp \p{-2c_b^2\min\{ C_bt, (C_bt)^\gamma \} (n_P\delta^{2\gamma+\frac{d}{\beta}} + n_Q\delta^{2+\frac{d}{\beta}})  }.
	\ealigns
	
	Note that
	\baligns
	n_P\delta^{2\gamma+\frac{d}{\beta}} + n_Q\delta^{2+\frac{d}{\beta}} &= \frac{n_P}{(n_P^{\frac{2\beta+d}{2\gamma\beta+d}}+n_Q)^{\frac{2\gamma\beta + d}{2\beta+d}}} + \frac{n_Q}{n_P^{\frac{2\beta+d}{2\gamma\beta+d}}+n_Q} \\
	&\geq \frac{n_P}{(n_P^{\frac{2\beta+d}{2\gamma\beta+d}}+n_P^{\frac{2\beta+d}{2\gamma\beta+d}})^{\frac{2\gamma\beta + d}{2\beta+d}}}\ind{n_Q \leq n_P^{\frac{2\beta+d}{2\gamma\beta+d}}} + \frac{n_Q}{n_Q+n_Q} \ind{n_Q > n_P^{\frac{2\beta+d}{2\gamma\beta+d}} } \\
	&= 2^{-\frac{2\gamma\beta+d}{2\beta+d}}\ind{n_Q \leq n_P^{\frac{2\beta+d}{2\gamma\beta+d}}} + \frac{1}{2}\ind{n_Q > n_P^{\frac{2\beta+d}{2\gamma\beta+d}} } \\
	&\geq \min\{ 2^{-\frac{2\gamma\beta+d}{2\beta+d}}, \frac{1}{2}  \}.
	\ealigns
	
	Combine two inequalities above, because $\min\{ C_bt, (C_bt)^\gamma \} = C_b^{1\wedge\gamma}t^{1\wedge\gamma}$, with constant $c_1 = 2c_b^2C_b^{1\wedge\gamma}\min\{ 2^{-\frac{2\gamma\beta+d}{2\beta+d}}, \frac{1}{2}  \}$,  we have under event $E_P\cap E_Q$, for all $x\in\Omega$,
	\beq \P_{Y|X}(\hf_{NN}(x) \neq f^\ast_Q(x)) \leq \exp \p{-c_1t^{1\wedge\gamma}}. \label{mcb2}  \eeq
	
	From lemma \ref{distb} we know
	\baligns
		\P(E_P \cap E_Q) &\geq 1-\P(E_P^c) - \P(E_Q^c) \\
		&\geq 1 - \frac{n_P}{k_P}\exp(-\frac{k_P}{6}) - \frac{n_Q}{k_Q}\exp(-\frac{k_Q}{6}) \\
		&\geq 1 - 2\delta^{1+\alpha}.
	\ealigns
	
	Therefore inequality $\eqref{mcb2}$ holds with probability larger than $1-2\delta^{1+\alpha}$.
	
	\ 
	
	\textbf{Case 2:} $\frac{n_Q}{k_Q}\exp(-\frac{k_Q}{6}) < \delta^{1+\alpha}$ and $\frac{n_P}{k_P}\exp(-\frac{k_P}{6}) \geq \delta^{1+\alpha}$
	
	Note that $\frac{k_P}{n_P} = \delta^{\frac{d}{\beta}}$, thus $\frac{n_P}{k_P}\exp(-\frac{k_P}{6}) \geq \delta^{1+\alpha}$ implies
	\[ k_P \leq 6(1+\alpha+\frac{d}{\beta})(-\log \delta).  \]
	
	And by $n_P = k_P\delta^{-\frac{d}{\beta}}$ we have
	\beq n_P \leq \frac{6\beta}{2\beta+d}(1+\alpha+\frac{d}{\beta})(n_P^{\frac{2\beta + d}{2\gamma\beta+d}} + n_Q)^\frac{d}{2\beta+d}\log(n_P^{\frac{2\beta + d}{2\gamma\beta+d}} + n_Q). \label{misc3} \eeq
	
	Then we further divide our discussion into two cases:
	
	\begin{enumerate}
		\item If $n_P^{\frac{2\beta + d}{2\gamma\beta+d}} \geq n_Q$ 
		
		Then \eqref{misc3} implies
		\baligns
			n_P &\leq \frac{6\beta}{2\beta+d}(1+\alpha+\frac{d}{\beta})(n_P^{\frac{2\beta + d}{2\gamma\beta+d}} + n_P^{\frac{2\beta + d}{2\gamma\beta+d}})^\frac{d}{2\beta+d}\log(n_P^{\frac{2\beta + d}{2\gamma\beta+d}} + n_P^{\frac{2\beta + d}{2\gamma\beta+d}}) \\
			&\leq \frac{6\beta}{2\beta+d}(1+\alpha+\frac{d}{\beta})2^\frac{d}{2\beta+d}n_P^{\frac{d}{2\gamma\beta+d}}2\frac{2\beta + d}{2\gamma\beta+d}\log n_P.
		\ealigns
		
		Note that $\frac{d}{2\gamma\beta+d} < 1$, so the above inequality implies $n_P$ cannot be larger than a certain constant depending on $\alpha,\beta,\gamma$ and $d$. Also by $n_Q < n_P^{\frac{2\beta + d}{2\gamma\beta+d}}$ we know $n_Q$ is bounded. Given $n_P$ and $n_Q$ are both bounded the inequality $\eqref{c_1}$ is trivial with large enough $C_1$.
		
		\item If $n_P^{\frac{2\beta + d}{2\gamma\beta+d}} < n_Q$
		
		We have
		\[ w_Pk_P = \delta^\gamma k_P \leq 6(1+\alpha+\frac{d}{\beta})(-\log \delta)\delta^\gamma < \frac{c_b}{4}. \]
		when $\delta$ is smaller than some constant depending on $\alpha,\beta,\gamma,d$ and $c_b$. (When $\delta$ is larger than this constant, the inequality $\eqref{c_1}$ is trivial with large enough $C_1$.)
		
		Also we have
		\[  w_Qk_Q\zeta_Q(x) \geq C_b\delta^{\frac{2\beta+d}{\beta}} n_Q \geq \frac{n_Q}{n_P^{\frac{2\beta + d}{2\gamma\beta+d}} + n_Q} \geq \frac{1}{2}. \]
		given $C_b > 1$ is large enough.
		
		Note that $\zeta_Q(x) < \frac{1}{2}$, given $c_b<1$ is small enough, combine two inequalities above we have
		\[ c_bw_Qk_Q\zeta_Q(x) - w_Pk_P \geq \frac{1}{3}(c_bw_Qk_Q\zeta_Q(x) + w_Pk_P) \geq \frac{c_b}{3}(w_Qk_Q\zeta_Q(x) + w_Pk_P\zeta_Q(x)^\gamma). \]
		
		Then based on \eqref{mcb1} Lemma $\ref{tailb}$, we have under event $E_Q$, given $C_b > 1$ large enough, for all $x$ that satisfies $\zeta_Q(x) \geq C_b\delta$, 
		\baligns \P_{Y|X}(\hf_{NN}(x) \neq f^\ast_Q(x)) &\leq \exp\p{ -2c_b^2\frac{[(w_Qk_Q\zeta_Q(x) - w_Pk_P)\vee 0]^2}{k_Pw_P^2 + k_Qw_Q^2} } \\
		&\leq \exp\p{ -\frac{2}{9}c_b^2\frac{(w_Qk_Q\zeta_Q(x) + w_Pk_P\zeta_Q(x)^\gamma)^2}{k_Pw_P^2 + k_Qw_Q^2} }.
		\ealigns
		
		Then follow similar proofs in case 1 we have with proper choice of $c_1$, inequality \eqref{mcb2} holds with probability at least
		\[ \P(E_Q) \geq 1 - \frac{n_Q}{k_Q}\exp(-\frac{k_Q}{6}) \geq 1 - \delta^{1+\alpha}.  \]
	\end{enumerate}
	\ 
	
	\textbf{Case 3:} $\frac{n_Q}{k_Q}\exp(-\frac{k_Q}{6}) \geq \delta^{1+\alpha}$ and $\frac{n_P}{k_P}\exp(-\frac{k_P}{6}) < \delta^{1+\alpha}$
	
	The proof is symmetric to the proof of case 2. Thus we omit the proof here.
	
	\ 
	
	\textbf{Case 4:} $\frac{n_Q}{k_Q}\exp(-\frac{k_Q}{6}) \geq \delta^{1+\alpha}$ and $\frac{n_P}{k_P}\exp(-\frac{k_P}{6}) \geq \delta^{1+\alpha}$
	
	In this case we still have inequality \eqref{misc3}
	
	\[ n_P \leq \frac{6\beta}{2\beta+d}(1+\alpha+\frac{d}{\beta})(n_P^{\frac{2\beta + d}{2\gamma\beta+d}} + n_Q)^\frac{d}{2\beta+d}\log(n_P^{\frac{2\beta + d}{2\gamma\beta+d}} + n_Q). \]
	
	And similarly $\frac{n_Q}{k_Q}\exp(-\frac{k_Q}{6}) \geq \delta^{1+\alpha}$ implies
	
	\[ n_Q \leq \frac{6\beta}{2\beta+d}(1+\alpha+\frac{d}{\beta})(n_P^{\frac{2\beta + d}{2\gamma\beta+d}} + n_Q)^\frac{d}{2\beta+d}\log(n_P^{\frac{2\beta + d}{2\gamma\beta+d}} + n_Q). \]
	
	Combine two inequalities above we have
	\[ n_P\vee n_Q \leq \frac{6\beta}{2\beta+d}(1+\alpha+\frac{d}{\beta})((n_P\vee n_Q)^{\frac{2\beta + d}{2\gamma\beta+d}} + (n_P\vee n_Q))^\frac{d}{2\beta+d}\log((n_P\vee n_Q)^{\frac{2\beta + d}{2\gamma\beta+d}} + (n_P\vee n_Q)).  \]
	
	which implies $n_P\vee n_Q$ is upper bounded by some constant. Therefore the inequality \eqref{c_1} is trivial with large enough $C_1$. 
	
	\qed

	\section{Proof of Theorem \ref{lb}}
	
	The proof of lower bound theorem \ref{lb} are mainly based on a construction of two families of distributions $P_\sigma$ and $Q_\sigma, \sigma \in \{0,1\}^m$ and applying Assouad's lemma on family $P_\sigma^{\otimes n_P}\times Q_\sigma^{\otimes n_Q}, \sigma \in \{0,1\}^m$. 
	
	First, let's define several quantities which will be used later in the proof. Let
	\[ r = c_r (n_P^{\frac{2\beta + d}{2\gamma\beta+d}} + n_Q)^{-\frac{1}{2\beta+d}}, \quad w = c_w r^d,\quad m = \lfloor c_m r^{\alpha\beta-d} \rfloor. \]
	where $c_r, c_w, c_m$ are universal constants which will be specified later. It is worthwhile to mention that as $n = n_P+n_Q \to \infty$, there will be $r, w \to 0$ and $m \to \infty$ because $\alpha\beta < d$.
	
	Let
	\[ G = \{ (6k_1r, 6k_2r, ..., 6k_dr), k_i = 1,2,..., \lfloor (6r)^{-1} \rfloor, i=1,2,...,d  \}. \]
	be a grid of $|G| = M = \lfloor (6r)^{-1} \rfloor^d$ points in the unit cube $\Omega$. Denote $x_1, x_2, ..., x_M$ as points in $G$.
	
	We are interested in $B(x_k, 2r), k=1,2,...,m$ balls with radius $2r$ centered at $x_k$. Let $B_c = \Omega \backslash \bigcup_{k=1}^m B(x_k,2r)$ denote the points that aren't in any of the $m$ balls. Note that $B(x_k, 2r), k=1,2,...,m$ are mutually disjoint, so $B(x_k, 2r), k=1,2,...,m$ and $B_c$ forms a partition of $\Omega$. A side note of above arguments: extracting $m$ center points out of $M$ is feasible because
	\[ m \approx c_m r^{\alpha\beta-d} < c_m r^{-d} < (6r)^{-d} \approx M  \]
	provided $c_m$ is small enough.

	Define function $g(\cdot)$ on $[0,\infty)$:
	\[ g(z) = \begin{cases} 1 \quad &0\leq z < 1 \\ 2-z \quad &1\leq z < 2 \\ 0 \quad &z\geq 2 \end{cases}. \]
	
	And define
	\[ h_Q(z) = C_\beta r^\beta g^\beta(z/r) \]
	\[ h_P(z) = C_\gamma C_\beta^\gamma r^{\beta\gamma}g^{\beta\gamma}(z/r).  \]
	
	By $r \leq c_r$ we have $\max(h_Q(z), h_P(z)) \leq \max(C_\beta c_r^\beta, C_\gamma C_\beta^\gamma c_r^{\beta\gamma})$. We choose $c_r$ small enough so that $\max(h_Q(z), h_P(z)) < 1$.
	
	Define the hypercube $\mathcal{H}$ of pairs $(P_\sigma, Q_\sigma)$ by
	\[ \mathcal{H} = \{(P_\sigma, Q_\sigma), \sigma = (\sigma_1,\sigma_2,...,\sigma_m) \in \{ -1,1 \}^m \}\]
	where both $P_\sigma, Q_\sigma$ are probability distributions of $(X,Y)$ on $R^d\times\{0,1\}$. We will construct each $(P_\sigma, Q_\sigma) \in \mathcal{H}$ by specifying conditional distributions $P_{\sigma, Y|X}, Q_{\sigma, Y|X}$ and marginal distributions $P_{\sigma,X}, Q_{\sigma,X}$.
	
	\textsc{Construction of $P_{\sigma, Y|X}$ and $Q_{\sigma, Y|X}$}:
	
	It is equivalent to specify regression functions $\eta_{P, \sigma}(x)$ and $\eta_{Q,\sigma}(x)$, defined as follows:
	\[ \eta_{P, \sigma}(x) = \begin{cases} \frac{1}{2}\p{1+\sigma_k h_P(\Vert x - x_k \Vert)} \quad &\text{if } x\in B(x_k, 2r) \text{ for some }k = 1,2,...,m \\ \frac{1}{2} \quad &\text{otherwise (equivalently $x \in B_c$).} \end{cases} \]
	\[ \eta_{Q, \sigma}(x) = \begin{cases} \frac{1}{2}\p{1+\sigma_k h_Q(\Vert x - x_k \Vert)} \quad &\text{if } x\in B(x_k, 2r) \text{ for some }k = 1,2,...,m \\ \frac{1}{2} \quad &\text{otherwise (equivalently $x \in B_c$).} \end{cases} \]
	
	\textsc{Construction of $P_{\sigma, X}$ and $Q_{\sigma, X}$}:
	
	Let $P_{\sigma, X}, Q_{\sigma, X}, \sigma\in\{0,1\}^m$ all have the same marginal distribution on $X$, with density $\mu(x)$. Define $\mu(x)$ as follows:
	\[ \mu(x) = \begin{cases} \frac{w}{\lambda[B(x_k,r)]} \quad &\text{if } x\in B(x_k, r) \text{ for some }k = 1,2,...,m \\ \frac{1-mw}{\lambda[B_c]} \quad &\text{if } x\in B_c \\ 0 \quad &\text{otherwise.} \end{cases} \]
	
	It is easy to verify that $\mu(x)$ is a density function on $\Omega$.
	
	Given the construction, next we are going to verify that distribution pairs $(P_\sigma, Q_\sigma) \in \mathcal{H}$ satisfies our assumptions, i.e.
	\[ (P_\sigma, Q_\sigma) \in \Pi(\alpha,\beta,\gamma,\mu) \text{ for all } (P_\sigma, Q_\sigma) \in \mathcal{H}.  \]
	
	\textsc{Verify Margin Assumption ($\alpha$)}:
	For any $\sigma \in \{-1,1\}^m$ we have
	\baligns &\P_\sigma(0<|\eta_{\sigma,Q}(X) - \frac{1}{2}|\leq t)  \\
	= & m\P_\sigma(0< h_Q(\Vert X-x_1 \Vert )\leq 2t) \\
	= & m \int_{B(x_1,r)} \ind{0< h_Q(\Vert X-x_1 \Vert )\leq 2t}\frac{w}{\lambda[B(x_k,r)]}dx \\
	= & mw \ind{t \geq C_\beta r^\beta/2} \\
	= & c_mc_w r^{\alpha\beta} \ind{t \geq C_\beta r^\beta/2} \leq C_\alpha t^\alpha.
	\ealigns
	
	provided that $c_m < C_\alpha (C_\beta/2)^\alpha/c_w$ is small enough. So $Q \in \mathcal{M}(\alpha, C_\alpha)$
	
	\textsc{Verify H\"{o}lder Smoothness ($\beta$)}: Note that $g(z)$ is 1-Lipschitz. For any $x, x' \in B(x_k, 2r)$, the basic inequality $|a^\beta - b^\beta| < |a-b|^\beta$ implies
	\baligns 
	\big|h_Q(\Vert x-x_k \Vert) - h_Q(\Vert x'-x_k \Vert)\big| &= C_\beta r^\beta \big|g^\beta(\Vert x-x_k \Vert/r) - g^\beta(\Vert x'-x_k \Vert/r)\big| \\
	&\leq  C_\beta r^\beta \big|g(\Vert x-x_k \Vert/r) - g(\Vert x'-x_k \Vert/r)\big|^\beta \\
	&\leq  C_\beta r^\beta \big|\Vert x-x_k \Vert/r - \Vert x'-x_k \Vert/r\big|^\beta \\
	&\leq C_\beta \Vert x - x' \Vert^\beta.
	\ealigns
	
	So $h_Q(\Vert x - x_k \Vert) \in \mathcal{H}(\beta, C_\beta)$. Tt is easy to extend to show that $\eta_{Q,\sigma}(x)\in \mathcal{H}(\beta, C_\beta)$.
	
	\textsc{Verify Relative Signal Exponent ($\gamma$)}: If $x\in B(x_k,2r)$ for some $k=1,2,...,m$, then $\sigma_k = 1$ suggests $\eta_{P, \sigma}(x) - \frac{1}{2} \geq 0$ and $\eta_{Q, \sigma}(x) - \frac{1}{2} \geq 0$; $\sigma_k = -1$ suggests $\eta_{P, \sigma}(x) - \frac{1}{2} \leq 0$ and $\eta_{Q, \sigma}(x) - \frac{1}{2} \leq 0$. If $x\in B_c$ then $\eta_{P, \sigma}(x) - \frac{1}{2} = \eta_{Q, \sigma}(x) - \frac{1}{2} = 0$. Therefore $\eqref{rse1}$ is verified.
	
	Also note that if $x\in B(x_k,2r)$ for some $k=1,2,...,m$, 
	\[ |\eta_{P, \sigma}(x) - \frac{1}{2}| = C_\gamma C_\beta^\gamma r^{\beta\gamma}g^{\beta\gamma}(\Vert x - x_k \Vert/r) = C_\gamma|\eta_{Q, \sigma}(x) - \frac{1}{2}|.  \]
	
	And if $x\in B_c$ we have
	\[ |\eta_{P, \sigma}(x) - \frac{1}{2}| = 0 = C_\gamma|\eta_{Q, \sigma}(x) - \frac{1}{2}|.  \]
	
	Therefore $\eqref{rse2}$ is verified. So $(P, Q) \in \Gamma(\gamma,C_\gamma)$

	\textsc{Verify Strong Density Assumption ($\mu$)}:
	It is easy to see that the support of $\mu_\sigma(x)$: $B_c \bigcup \p{\bigcup_{k=1}^m B(x_k, r)}$ is regular.
	
	If $x \in B_c$ we have
	\[ \mu(x) = \frac{1-mw}{1-m\lambda[B(x_1, 2r)]} = \frac{1-c_mc_wr^{\alpha\beta}}{1-2^d c_m \pi_d r^{\alpha\beta}} = 1+o(1).  \]
	
	If $x\in B(x_k, r) \text{ for some }k = 1,2,...,m$ we have
	\[ \mu(x) = \frac{w}{\lambda[B(x_1,r)]} = \pi_d^{-1}c_w.  \]
	
	Thus the marginal distribution satisfies strong density assumption with $\mu$ provided that $\pi_d\mu_- < c_w < \pi_d \mu_+$.
	
	So now we can conclude that the hypercube $\mathcal{H} \subset \Pi(\alpha,\beta,\gamma,\mu)$ with proper choices of $c_r, c_m$ and $c_w$.
	
	Finally, we are going to apply Assouad's lemma to proof the lower bound. let $H(\cdot, \cdot)$ denote the Hellinger distance between two measures. If $\sigma, \sigma' \in \{0,1\}^n$ are two indices that differ only at one element, i.e. $\sigma_k \neq \sigma_k'$ for some $k$ and $\sigma_i = \sigma_i'$ for all $i\neq k$. We have
	\baligns H^2(Q_\sigma, Q_{\sigma'}) &= \frac{1}{2}\int \mu(x) \p{ \p{\sqrt{\eta_{Q,\sigma}(x)} - \sqrt{\eta_{Q,\sigma'}(x)}}^2  + \p{\sqrt{1-\eta_{Q,\sigma}(x)} - \sqrt{1-\eta_{Q,\sigma'}(x)}}^2 }dx \\
	&= \frac{1}{2}\int_{B(x_k, r)} \frac{w}{\lambda[B(x_k,r)]}\cdot 2\p{\sqrt{\frac{1}{2}(1+C_\beta r^\beta)} - \sqrt{\frac{1}{2}(1-C_\beta r^\beta)}}^2 dx \\
	&= w(1 - \sqrt{1 - C_\beta^2 r^{2\beta}}) \\
	&\leq C_\beta^2 wr^{2\beta}.
	\ealigns
	
	Similarly we have
	\baligns H^2(P_\sigma, P_{\sigma'}) &= \frac{1}{2}\int \mu(x) \p{ \p{\sqrt{\eta_{P,\sigma}(x)} - \sqrt{\eta_{P,\sigma'}(x)}}^2 + \p{\sqrt{1-\eta_{P,\sigma}(x)} - \sqrt{1-\eta_{P,\sigma'}(x)}}^2 }dx \\
	&= \frac{1}{2}\int_{B(x_k, r)} \frac{w}{\lambda[B(x_k,r)]}\cdot 2\p{\sqrt{\frac{1}{2}(1+C_\gamma C_\beta^{\gamma} r^{\beta\gamma})} - \sqrt{\frac{1}{2}(1-C_\gamma C_\beta^{\gamma} r^{\beta\gamma})}}^2 dx \\
	&= w(1 - \sqrt{1 - C_\gamma^2 C_\beta^{2\gamma} r^{2\beta\gamma}}) \\
	&\leq C_\gamma^2 C_\beta^{2\gamma} w r^{2\beta\gamma}.
	\ealigns
	
	Therefore we have
	\baligns
		H^2(P_\sigma^{\otimes n_P}\times Q_\sigma^{\otimes n_Q}, &P_{\sigma'}^{\otimes n_P}\times Q_{\sigma'}^{\otimes n_Q} ) \leq n_P H^2(P_\sigma, P_{\sigma'}) + n_Q H^2(Q_\sigma, Q_{\sigma'}) \\
		\leq& C_\gamma^2 C_\beta^{2\gamma} c_w n_Pr^{2\beta\gamma+d} + C_\beta^2 c_w n_Q r^{2\beta+d} \\
		\leq& C_\gamma^2 C_\beta^{2\gamma} c_w c_r^{2\gamma\beta+d} n_P(n_P^{\frac{2\beta + d}{2\gamma\beta+d}} + n_Q)^{-\frac{2\beta\gamma+d}{2\beta+d}} + C_\beta^2c_wc_r^{2\beta+d}n_Q(n_P^{\frac{2\beta + d}{2\gamma\beta+d}} + n_Q)^{-1} \\
		\leq& \max(C_\gamma^2 C_\beta^{2\gamma} c_w c_r^{2\gamma\beta+d}, C_\beta^2c_wc_r^{2\beta+d}) \\& \cdot\p{n_P(n_P^{\frac{2\beta + d}{2\gamma\beta+d}} + n_Q)^{-\frac{2\beta\gamma+d}{2\beta+d}} + n_Q(n_P^{\frac{2\beta + d}{2\gamma\beta+d}} + n_Q)^{-1}} \\
		\leq& 2\max(C_\gamma^2 C_\beta^{2\gamma} c_w c_r^{2\gamma\beta+d}, C_\beta^2c_wc_r^{2\beta+d}) \\
		\leq& \frac{1}{4}.
	\ealigns
	provided that $c_r$ is small enough (doesn't depend on choice of $\sigma, \sigma'$).
	
	The above bound on hellinger distance implies that 
	\beq TV(P_\sigma^{\otimes n_P}\times Q_\sigma^{\otimes n_Q}, P_{\sigma'}^{\otimes n_P}\times Q_{\sigma'}^{\otimes n_Q} ) \leq \sqrt{2}H(P_\sigma^{\otimes n_P}\times Q_\sigma^{\otimes n_Q}, P_{\sigma'}^{\otimes n_P}\times Q_{\sigma'}^{\otimes n_Q} ) \leq \frac{\sqrt{2}}{2}. \label{assouad1} \eeq
	
	Also, for any classfier $\hat f$ we have
	\begin{equation}\begin{aligned} R_{Q_\sigma}(\hat f) + R_{Q_{\sigma'}}(\hat f) =& 2\E_{X \sim Q_{X,\sigma}}(|\eta_Q(X) - \frac{1}{2}|\ind{\hf(X) = f^\ast_{Q_\sigma}(X)})  \\ &+ 2\E_{X \sim Q_{X,\sigma'}}(|\eta_Q(X) - \frac{1}{2}|\ind{\hf(X) = f^\ast_{Q_{\sigma'}}(X)}) \\
	=& 2\sum_{i=1}^m \int_{B(x_i, r)} \mu(x) \cdot \frac{1}{2}C_\beta r^\beta \cdot \p{\ind{\hf(x) = f^\ast_{Q_\sigma}(x)} + \ind{\hf(x) = f^\ast_{Q_{\sigma'}}(x)}}dx \\
	\geq& \int_{B(x_k, r)} \mu(x) \cdot C_\beta r^\beta \cdot \p{\ind{\hf(x) = f^\ast_{Q_\sigma}(x)} + \ind{\hf(x) = f^\ast_{Q_{\sigma'}}(x)}}dx \\
	=& C_\beta wr^\beta. \label{assouad2}
	\end{aligned}\end{equation}
	Because when $x \in B(x_k, r)$ we have $f^\ast_{Q_\sigma}(x) \neq f^\ast_{Q_{\sigma'}}(x)$ thus $\ind{\hf(x) = f^\ast_{Q_\sigma}(x)} + \ind{\hf(x) = f^\ast_{Q_{\sigma'}}(x)} = 1$.
	
	Finally, apply the Assouad's lemma based on established inequalities $\eqref{assouad1}$ and $\eqref{assouad2}$, we have for all estimators $\hat f$,
	\[ \max_{(P,Q)\in \mathcal{H}} \mathcal{E}_Q(\hat f) \geq \frac{m}{2}\cdot C_\beta wr^\beta \cdot (1-\frac{\sqrt{2}}{2}) = \frac{2-\sqrt{2}}{4}C_\beta c_m c_w c_r^{(1+\alpha)\beta}(n_P^{\frac{2\beta + d}{2\gamma\beta+d}} + n_Q)^{-\frac{\beta(1+\alpha)}{2\beta+d}}. \]
	which gives the minimax lower bound. \qed

	\section{Proof of Theorem \ref{adap_bound}}
	
	First we give an auxillary lemma showing an union bound on difference between any weighted $K$-NN estimator and its mean.
	
	\begin{lemma} \label{adapt_lemma_1}
		Define weighted $K$-NN estimator $\hat\eta_{k_P, k_Q ,w}(x)$ as
		\[  \hat\eta_{k_P, k_Q ,w}(x) = w\frac{\sum_{i=1}^{k_P}Y^P_{(i)}(x)}{k_P} + (1-w)\frac{\sum_{i=1}^{k_Q}Y^Q_{(i)}(x)}{k_Q}.  \]
		
		Then with probability at least $1-\delta$, for all $w\in [0,1], k_P\in [n_P], k_Q\in [n_Q], x\in \reals^d$, we have
		\[ | \hat\eta_{k_P, k_Q,w}(x) - \E_{Y|X}\hat\eta_{k_P, k_Q,w}(x) | \leq \sqrt{\bigg((d+1)\log(n_P+n_Q) - \log(\delta/2)\bigg)\p{\frac{w^2}{k_P} + \frac{(1-w)^2}{k_Q}}}.  \]
		
		Take $\delta = 2(n_P+n_Q)^{-2}$ we have with probability at least $1 - 2(n_P+n_Q)^{-2}$, for all $w\in [0,1], k_P\in [n_P], k_Q\in [n_Q], x\in \reals^d$, 
		\beq | \hat\eta_{k_P, k_Q,w}(x) - \E_{Y|X}\hat\eta_{k_P, k_Q,w}(x) | \leq \sqrt{(d+3)\log(n_P+n_Q) \p{\frac{w^2}{k_P} + \frac{(1-w)^2}{k_Q}}}. \label{ada_lm} \eeq
		
	\end{lemma}

	The proof of Lemma $\ref{adapt_lemma_1}$ is provided in Section \ref{sec.lemma5}. In the following proofs, we define the event $E_A$ be that the inequality $\eqref{ada_lm}$ holds for all $w, k_P, k_Q$ stated in the lemma \ref{adapt_lemma_1}. Thanks to the above lemma we have
	\[ \P(E_A) \geq 1 - 2(n_P+n_Q)^{-2}. \]
	
	First of all, let's define some important quantities. Let
	\[ \delta = C_\delta \p{\p{\frac{n_P}{\log (n_P + n_Q)}}^{\frac{2\beta + d}{2\gamma\beta+d}} + \frac{n_Q}{\log (n_P + n_Q)}}^{-\frac{\beta}{2\beta+d}}  \]
	where $C_\delta > 0$ is a large constant which will be given later. Define $G_\delta$ be the set
	\[ G_\delta = \{ x : \zeta_Q(x) \geq \delta \}.  \]
	
	Also we define $k^{opt}(x), k_P^{opt}(x), k_Q^{opt}(x)$ be the "optimal" (oracle) choices of number of neighbors defined by
	\[ k^{opt}(x) = \max_{\Vert X_{(k)}(x) - x \Vert \leq ( \frac{\delta}{2C_\beta} )^{1/\beta} } k, \]
	\[ k_P^{opt}(x) = k_P^{(k^{opt}(x))},  \quad k_Q^{opt}(x) = k_Q^{(k^{opt}(x))}.  \]
	
	where $k_P^{(k^{opt}(x))}$ ($k_Q^{(k^{opt}(x))}$) is number of covariates from $P-$data ($Q-$data) among all $k^{opt}(x)$ nearest covariates to $x$, as is defined in Algorithm \ref{ada_algo}. We will sometimes omit $x$ and just write $k^{opt}, k_P^{opt}, k_Q^{opt}$ if no confusion in the context.
	
	Define $k^{stop}(x)$ be the stopping time of the algorithm \ref{ada_algo}. i.e.
	\[ k^{stop}(x) = \min_{\hat r^{(k)} > (d+3)\log(n_P+n_Q)} k.  \]
	
	And let $k^{stop}(x) = \infty$ if the algorithm doen't stop till the end. Similarly denote
	\[ k_P^{stop}(x) = k_P^{(k^{stop}(x))}, \quad k_Q^{stop}(x) = k_Q^{(k^{stop}(x))}.  \]
	
	And sometimes we will omit $x$ for simplicity.
	
	Next we are going to state several claims leading to prove theorem \ref{adap_bound} step by step. After each claim we will directly provide a proof of that claim. 
		
	\begin{claim}
		If $x \in G_\delta$, $k_P \leq k_P^{opt}(x)$, $k_Q \leq k_Q^{opt}(x)$, we have
		\begin{enumerate}
			\item When $f^\ast(x) = 1$,
			\beq \E_{Y|X}(\bar Y^P_{(1:k_P)}(x)) - \frac{1}{2} \geq c_b\delta^\gamma \label{cl1_1}  \eeq
			\beq \E_{Y|X}(\bar Y^Q_{(1:k_Q)}(x)) - \frac{1}{2} \geq c_b\delta \label{cl1_2}.  \eeq
			\item When $f^\ast(x) = 0$,
			\beq \E_{Y|X}(\bar Y^P_{(1:k_P)}(x)) - \frac{1}{2} \leq -c_b\delta^\gamma \label{cl1_3} \eeq
			\beq \E_{Y|X}(\bar Y^Q_{(1:k_Q)}(x)) - \frac{1}{2} \leq -c_b\delta \label{cl1_4}. \eeq
		\end{enumerate}
		\label{claim1}
	\end{claim}

	\textsc{Proof of claim \ref{claim1}}. Note that if $x \in G_\delta$, $k_P \leq k_P^{opt}$, we have
	\[ \Vert X^P_{(k_P)}(x) - x \Vert \leq \Vert X^P_{(k_P^{opt})}(x) - x \Vert \leq \Vert X_{(k^{opt})}(x) - x \Vert \leq ( \frac{\delta}{2C_\beta} )^{1/\beta}.  \]
	
	Therefore we have
	\[ \zeta_Q(x) \geq \delta \geq 2C_\beta \Vert X^P_{(k_P)}(x) - x \Vert^\beta. \]
	
	Apply lemma \ref{biasb} we can obtain \eqref{cl1_1} and \eqref{cl1_3} given $\zeta_Q(x) \geq \delta > 0$
	
	Similarly, note that if $x \in G_\delta$, $k_Q \leq k_Q^{opt}$, we have
	\[ \Vert X^Q_{(k_Q)}(x) - x \Vert \leq \Vert X^Q_{(k_Q^{opt})}(x) - x \Vert \leq \Vert X_{(k^{opt})}(x) - x \Vert \leq ( \frac{\delta}{2C_\beta} )^{1/\beta}.  \]
	
	Therefore we have
	\[ \zeta_Q(x) \geq \delta \geq 2C_\beta \Vert X^Q_{(k_Q)}(x) - x \Vert^\beta. \]

	Apply lemma \ref{biasb} we can obtain \eqref{cl1_2} and \eqref{cl1_4} given $\zeta_Q(x) \geq \delta > 0$.
	
	\qed
	
	\begin{claim}
		Under event $E_A$, if $x \in G_\delta$ and $k^{stop}(x) \leq k^{opt}(x)$, then the output of algorithm is correct, i.e. $\hat f_a(x) = f^\ast(x)$.
		\label{claim2}
	\end{claim}

	\textsc{Proof of claim \ref{claim2}}. $k^{stop}(x) \leq k^{opt}(x)$ implies $k^{stop} < \infty$ so the algorithm stops at the round $k^{stop}$. By the stopping rule we know that
	\[ \hat r^{(k^{stop})} > (d+3)\log(n_P+n_Q).  \]
	
	By construction of $\hat r^{(k^{stop})}$ it is easy to show that
	\[ \sqrt{r^{(k^{stop})}} = \max_{w} \frac{\left| w (\bar Y^P_{(1:k_P^{stop})}(x) - \frac{1}{2}) + (1-w) (\bar Y^Q_{(1:k_Q^{stop})}(x) - \frac{1}{2}) \right|}{\sqrt{\frac{w^2}{k_P^{stop}} + \frac{(1-w)^2}{k_Q^{stop}}}}.  \]
	
	Let $w_0$ be one of the value of $w$ such that right hand side takes its maximum. Combine two formulas above we have
	\[ \left| w_0 \bar Y^P_{(1:k_P^{stop})}(x) + (1-w_0) \bar Y^Q_{(1:k_Q^{stop})}(x) - \frac{1}{2} \right| > \sqrt{(d+3)\log(n_P+n_Q)(\frac{w_0^2}{k_P^{stop}} + \frac{(1-w_0)^2}{k_Q^{stop}})}. \]
	
	For simplicity we may rewrite the left hand side as $|\hat\eta_{k_P^{stop}, k_Q^{stop} ,w_0}|$ as is defined in lemma \ref{adapt_lemma_1}. By definition of $E_A$ we know that under $E_A$ we have
	\[ |\hat\eta_{k_P^{stop}, k_Q^{stop} ,w_0}(x) - \E_{Y|X}\hat\eta_{k_P^{stop}, k_Q^{stop} ,w_0}(x)| \leq \sqrt{(d+3)\log(n_P+n_Q)(\frac{w_0^2}{k_P^{stop}} + \frac{(1-w_0)^2}{k_Q^{stop}})}.  \]
	
	Combine two inequalities above we have
	\[ |\hat\eta_{k_P^{stop}, k_Q^{stop} ,w_0}(x) - \frac{1}{2}| > |(\hat\eta_{k_P^{stop}, k_Q^{stop} ,w_0}(x) - \frac{1}{2}) - (\E_{Y|X}\hat\eta_{k_P^{stop}, k_Q^{stop} ,w_0}(x)- \frac{1}{2})|, \]
	which implies
	\[ \sign(\hat\eta_{k_P^{stop}, k_Q^{stop} ,w_0}(x) - \frac{1}{2}) = \sign(\E_{Y|X}\hat\eta_{k_P^{stop}, k_Q^{stop} ,w_0}(x)- \frac{1}{2}) \neq 0.  \]
	
	Note that $k^{stop}(x) \leq k^{opt}(x)$ implies $k_P^{stop}(x) \leq k_Q^{opt}(x)$ and $k_Q^{stop}(x) \leq k_Q^{opt}(x)$, given $x \in G_\delta$, by claim \ref{claim1} we have when $f^\ast(x) = 1$,
	\[ \E_{Y|X}\hat\eta_{k_P^{stop}, k_Q^{stop} ,w_0}(x)- \frac{1}{2} \geq c_b(w_0\delta^\gamma + (1-w_0)\delta) > 0.  \]
	
	And when $f^\ast(x) = 0$, we have
	\[ \E_{Y|X}\hat\eta_{k_P^{stop}, k_Q^{stop} ,w_0}(x)- \frac{1}{2} \leq -c_b(w_0\delta^\gamma + (1-w_0)\delta) < 0.  \]
	
	So,
	\[ \sign(\hat\eta_{k_P^{stop}, k_Q^{stop} ,w_0}(x) - \frac{1}{2}) = \sign(\E_{Y|X}\hat\eta_{k_P^{stop}, k_Q^{stop} ,w_0}(x)- \frac{1}{2}) = \begin{cases} 1 \eqif f^\ast(x) = 1\\ -1 \eqif f^\ast(x) = 0 \end{cases}. \]
	
	By simple calculation one can show that
	\[ \sign(\hat\eta_{k_P^{stop}, k_Q^{stop} ,w_0}(x) - \frac{1}{2}) = \sign(\sqrt{k_P^{stop}}\p{Y^P_{(1:k_P^{stop})}(x) - \frac{1}{2}} + \sqrt{k_Q^{stop}} \p{Y^Q_{(1:k_Q^{stop})}(x)- \frac{1}{2}}) \neq 0. \]
	
	Therefore,
	\[ \hat f_a(x) = \ind{\sqrt{k_P^{stop}}(Y^P_{(1:k_P^{stop})}(x) - \frac{1}{2}) + \sqrt{k_Q^{stop}} (Y^Q_{(1:k_Q^{stop})}(x)- \frac{1}{2}) \geq 0} = \ind{\hat\eta_{k_P^{stop}, k_Q^{stop} ,w_0}(x) - \frac{1}{2} \geq 0} = f^\ast(x).  \]
	
	\qed
	
	\begin{claim}
		There exist a constant $C_2 > 0$ such that, with probability at least $1-C_2\delta^{1+\alpha}$, we have $\hat f_a(x) = f^\ast(x)$ for all $x \in G_\delta$.
		\label{claim3}
	\end{claim}

	\textsc{Proof of claim \ref{claim3}}. Here we divide our discussion into two cases: 
	
	Case 1: $\p{\frac{n_P}{\log (n_P + n_Q)}}^{\frac{2\beta + d}{2\gamma\beta+d}} \leq \frac{n_Q}{\log (n_P + n_Q)}$
	
	Case 2: $\p{\frac{n_P}{\log (n_P + n_Q)}}^{\frac{2\beta + d}{2\gamma\beta+d}} > \frac{n_Q}{\log (n_P + n_Q)}$
	
	The proofs for the above two cases are symmetric so here we only discuss the first case. Thus until the end of the proof of claim \ref{claim3} we assume $\p{\frac{n_P}{\log (n_P + n_Q)}}^{\frac{2\beta + d}{2\gamma\beta+d}} \leq \frac{n_Q}{\log (n_P + n_Q)}$.
	
	Let $\underline{k}_Q = \frac{n_Q}{C_D^d}(\frac{\delta}{2C_\beta})^{\frac{d}{\beta}}$. Apply lemma \ref{distb} with $k_Q = \underline{k}_Q$, we know that with probability at least $1-C_D\frac{n_Q}{\underline{k}_Q}\exp(-\frac{\underline{k}_Q}{6})$, for all $x$ we have
	\[ \Vert X^Q_{\underline{k}_Q}(x) - x \Vert \leq C_D(\frac{\underline{k}_Q}{n_Q})^\frac{1}{d}.  \]
	
	We still let $E_Q$ denote the event that above inequality holds for all $x$. Now we have
	
	\[ \P(E_Q) \geq 1-C_D\frac{n_Q}{\underline{k}_Q}\exp(-\frac{\underline{k}_Q}{6}).  \]
	
	Under event $E_Q$, we have for all $x$,
	\[ \Vert X^Q_{\underline{k}_Q}(x) - x \Vert \leq C_D(\frac{\underline{k}_Q}{n_Q})^\frac{1}{d} = (\frac{\delta}{2C_\beta})^\frac{1}{\beta} < \Vert X^Q_{k_Q^{opt}+1}(x) - x \Vert,  \]
	which implies
	\[ k_Q^{opt}(x) \geq \underline{k}_Q = \frac{n_Q}{C_D^d}(\frac{\delta}{2C_\beta})^{\frac{d}{\beta}} \text{ for all } x. \]
	
	Therefore, under event $E_Q$, for all $x \in G_\delta \cap \{ f^\ast(x) = 1\}$, apply claim \ref{claim1} we have
	\begin{equation}
		\begin{aligned}
		\sqrt{k_Q^{opt}(x)}&\p{\E_{Y|X}(\bar Y^Q_{(1:k_Q^{opt})}) - \frac{1}{2}} \geq \sqrt{\frac{n_Q}{C_D^d}(\frac{\delta}{2C_\beta})^{\frac{d}{\beta}}} \cdot c_b \delta \\
		& = \sqrt{\frac{c_b^2C_\delta^{2+d/\beta}}{C_D^d(2C_\beta)^{d/\beta}}n_Q\p{\p{\frac{n_P}{\log (n_P + n_Q)}}^{\frac{2\beta + d}{2\gamma\beta+d}} + \frac{n_Q}{\log (n_P + n_Q)}}^{-1}} \\
		&\geq \sqrt{\frac{c_b^2C_\delta^{2+d/\beta}}{C_D^d(2C_\beta)^{d/\beta}}n_Q\p{2\frac{n_Q}{\log (n_P + n_Q)}}^{-1}} \\
		&\geq 3\sqrt{(d+3)\log(n_P+n_Q)}
		\end{aligned}
		\label{cl2_1}
	\end{equation}
	with large enough choice of constant $C_\delta$.
	
	In addition, under event $E_A$, with choice of $w = 0$ and $k_Q = k_Q^{opt}(x)$, we have for all $x$,
	\beq |\bar Y^Q_{(1:k_Q^{opt})} - \E_{Y|X}(\bar Y^Q_{(1:k_Q^{opt})})| \leq \sqrt{\frac{(d+3)\log(n_P+n_Q)}{k_Q^{opt}}}. \label{cl2_2} \eeq
	
	Combine \eqref{cl2_1} and \eqref{cl2_2} together, under $E_Q\cap E_A$, for all $x \in G_\delta \cap \{ f^\ast(x) = 1\}$ we have
	\baligns
		\sqrt{k_Q^{opt}(x)}\p{\bar Y^Q_{(1:k_Q^{opt})} - \frac{1}{2}} \geq& \sqrt{k_Q^{opt}(x)}\p{\E_{Y|X}(\bar Y^Q_{(1:k_Q^{opt})}) - \frac{1}{2}} \\& \quad - \sqrt{k_Q^{opt}(x)}|\bar Y^Q_{(1:k_Q^{opt})} - \E_{Y|X}(\bar Y^Q_{(1:k_Q^{opt})})| \\
		\geq& 3\sqrt{(d+3)\log(n_P+n_Q)} - \sqrt{(d+3)\log(n_P+n_Q)} \\
		>& \sqrt{(d+3)\log(n_P+n_Q)}.
	\ealigns
	
	Apply similar derivation, we can obtain that under $E_Q\cap E_A$, for all $x \in G_\delta \cap \{ f^\ast(x) = 0\}$ we have
	\[ \sqrt{k_Q^{opt}(x)}\p{\bar Y^Q_{(1:k_Q^{opt})} - \frac{1}{2}} \leq - \sqrt{(d+3)\log(n_P+n_Q)}.  \]
	
	Therefore, under $E_Q\cap E_A$, for all $x \in G_\delta$ we have
	\[ k_Q^{opt}(x)\p{\bar Y^Q_{(1:k_Q^{opt})} - \frac{1}{2}}^2 > (d+3)\log(n_P+n_Q).  \]
	
	Note that $\hat r^{(k^{opt})} > k_Q^{opt}(x)\p{\bar Y^Q_{(1:k_Q^{opt})} - \frac{1}{2}}^2$ so
	\[ \hat r^{(k^{opt})} > (d+3)\log(n_P+n_Q). \]
	
	This means that the algorithm \ref{ada_algo} must stop at the round $k = k^{opt}(x)$ if it does not stop earlier.
	
	Therefore, under $E_Q\cap E_A$, for all $x \in G_\delta$, we have
	\[ k^{stop}(x) \leq k^{opt}(x). \]
	
	Now apply claim \ref{claim2} we know $k^{stop}(x) \leq k^{opt}(x)$ implies
	\[ \hat f_a(x) = f^\ast(x). \]
	for all $x\in G_\delta$ under event $E_Q\cap E_A$.
	
	So it remains to show that probability of the above event is at least $1-C_2\delta^{1+\alpha}$, i.e. we are going to show there exist a constant $C_2 > 0$ such that
	\[ \P(E_Q\cap E_A) \geq 1-C_2\delta^{1+\alpha}. \]
	
	It suffices to show that there exist some constants $C_{21}, C_{22} > 0$ such that
	\[ \P(E_A^c) \leq C_{21}\delta^{1+\alpha} \quad  \P(E_Q^c) \leq C_{22}\delta^{1+\alpha}. \]
	
	First, by lemma \ref{adapt_lemma_1} we have
	\[ \P(E_A^c) \leq 2(n_P + n_Q)^{-2}. \]
	
	Note that $\alpha\beta \geq d$ we have
	\[ \max(\frac{2\beta+d}{2\gamma\beta+d}\cdot \frac{\beta(1+\alpha)}{2\beta+d}, \frac{\beta(1+\alpha)}{2\beta+d}) < \frac{\beta + \beta\alpha}{d} \leq \frac{1 + d}{d} \leq 2. \]
	
	So
	\[ \P(E_A^c) \leq 2(n_P + n_Q)^{-2} < 2(n_P^\frac{2\beta+d}{2\gamma\beta+d} + n_Q)^{-\frac{\beta(1+\alpha)}{2\beta+d}} < 2\delta^{1+\alpha}.  \]
	
	Then we are going to bound $\P(E_Q)$. By lemma \ref{distb} we have
	\[ \P(E_Q^c) \leq C_D\frac{n_Q}{\underline{k}_Q}\exp(-\frac{\underline{k}_Q}{6}). \]
	
	And by $\p{\frac{n_P}{\log (n_P + n_Q)}}^{\frac{2\beta + d}{2\gamma\beta+d}} \leq \frac{n_Q}{\log (n_P + n_Q)}$ we have
	\baligns
		n_Q > \frac{n_Q}{\log(n_P+n_Q)} \geq \frac{1}{4}\p{\frac{n_Q}{\log(n_P+n_Q)} + \p{\frac{n_P}{\log (n_P + n_Q)}}^{\frac{2\beta + d}{2\gamma\beta+d}}} = \frac{1}{4}\delta^{-(2+\frac{d}{\beta})}.
	\ealigns
	
	Thus
	\baligns 
		\frac{\P(E_Q^c)}{\delta^{1+\alpha}} &\leq C_D^{d+1}(2C_\beta)^\frac{d}{\beta}\delta^{-(1+\alpha+d/\beta)}\exp(-\frac{n_Q}{6C_D^d}(\frac{\delta}{2C_\beta})^{\frac{d}{\beta}}) \\
		&\leq C_D^{d+1}(2C_\beta)^\frac{d}{\beta}\delta^{-(1+\alpha+d/\beta)}\exp(-\frac{1}{24C_D^d(2C_\beta)^\frac{d}{\beta}}\delta^{-2}).
	\ealigns
	
	The right hand side goes to zero as $\delta \to 0$, so it is bounded by some large enough constant $C_{22}$. Thus we have
	\[ \P(E_Q^c) \leq C_{22}\delta^{1+\alpha}. \]
	
	Now the proof of claim \ref{claim3} is completed.
	
	\qed
	
	\textsc{Proof of theorem \ref{adap_bound}}. Let the event $E_W$ be that $\hat f_a(x) = f^\ast(x)$ for all $x \in G_\delta$. From claim \ref{claim3} we know 
	\[ \P(E_W^c) \leq C_2 \delta^{1+\alpha}. \]
	
	Note that under $E_Z$, $\hat f_a(x) \neq f^\ast(x)$ implies that $x \notin G_\delta$. So we have
	\baligns
		\E \mathcal{E}_Q(\hat f_a) &= \E [\E_{X\sim Q_X}(\zeta_Q(X)\ind{\hat f_a(X) \neq f^\ast(X)})] \\
		&\leq \E [\E_{X\sim Q_X}(\zeta_Q(X)\ind{\hat f_a(X) \neq f^\ast(X)}) | E_W] + 1\cdot \P(E_W^c) \\
		&\leq \E [\E_{X\sim Q_X}(\zeta_Q(X)\ind{X \notin G_\delta}) | E_W] + \P(E_W^c) \\
		&= \E_{X\sim Q_X}(\zeta_Q(X)\ind{\zeta_Q(X) < \delta}) + \P(E_W^c) \\
		&\leq \delta\P_{X\sim Q_X}(\zeta_Q(X) < \delta) + \P(E_W^c) \\
		&\leq \delta \cdot C_\alpha\delta^\alpha + C_2 \delta^{1+\alpha} \\
		& = C_\delta(C_\alpha + C_2) \p{\p{\frac{n_P}{\log (n_P + n_Q)}}^{\frac{2\beta + d}{2\gamma\beta+d}} + \frac{n_Q}{\log (n_P + n_Q)}}^{-\frac{\beta(1+\alpha)}{2\beta+d}}.
	\ealigns
	
	\qed

	\section{Proof of Lemma \ref{adapt_lemma_1}} \label{sec.lemma5}
	
	Note that conditional on $X_{1:n}$, $Y_{(i)}^P(x), i=1,2,...,n_P$ are independent Bernuolli variables. So by Hoeffding's inequality
	\beq \P_{Y|X}( | \bar Y^P_{(1:k_P)}(x) - \E_{Y|X}\bar Y^P_{(1:k_P)}(x) | > \epsilon) \leq 2\exp\p{-2k_P\epsilon^2}. \label{ada_lm1_1}  \eeq
	
	For the same reason about $Q-$data we also have
	\beq \P_{Y|X}( | \bar Y^Q_{(1:k_Q)}(x) - \E_{Y|X}\bar Y^Q_{(1:k_Q)}(x) | > \epsilon) \leq 2\exp\p{-2k_Q\epsilon^2}. \label{ada_lm1_2}  \eeq
	
	Note that for any $x$, $X^P_{(1)}(x), X^P_{(2)}(x), ..., X^P_{(k_P)}(x)$ form a set of points that falls in the ball $B(x, \Vert X^P_{(k_P)}(x) - x \Vert )$. It is well-known that total number of sets of form $A = \{X_1^P, X_2^P, ..., X_{n_P}^P\} \cap B(x, r)$, with $x\in \reals^d, r\geq0$, i.e. the number of ways $Q-$covariates intercept with a ball, is upper-bounded by $n_Q^{d+1}$. This implies there are at most $n_P^{d+1}$ possible different random variables of form $\bar Y^P_{(1:k_P)}(x)$ with $x \in \reals^d$ and $k_P \in [n_P]$. For Similar reason, there are at most $n_Q^{d+1}$ possiblities of $\bar Y^Q_{(1:k_Q)}(x)$. 
	
	Plug $\epsilon = \sqrt{\frac{(d+1)\log(n_P+n_Q) - \log(\delta/2)}{2k_P}}$ into \eqref{ada_lm1_1}, plug $\epsilon = \sqrt{\frac{(d+1)\log(n_P+n_Q) - \log(\delta/2)}{2k_Q}}$ into \eqref{ada_lm1_2}, and apply the union bound, we have
	\baligns \P_{Y|X} \bigg(\exists x, k_P, k_Q \text{ s.t. }& | \bar Y^P_{(1:k_P)}(x) - \E_{Y|X}\bar Y^P_{(1:k_P)}(x) | > \sqrt{\frac{(d+1)\log(n_P+n_Q) - \log(\delta/2)}{2k_P}} \\ 
	\text{or }&  \bar Y^Q_{(1:k_Q)}(x) - \E_{Y|X}\bar Y^Q_{(1:k_Q)}(x) | > \sqrt{\frac{(d+1)\log(n_P+n_Q) - \log(\delta/2)}{2k_Q}}\bigg) \\
	\leq n_P^{d+1} \cdot \frac{\delta}{(n_P+n_Q)^{d+1}} + &n_Q^{d+1} \cdot \frac{\delta}{(n_P+n_Q)^{d+1}} \leq \delta.
	\ealigns
	
	So with probability at least $1-\delta$ we have for all $x, k_P, k_Q$
	\beq |\bar Y^P_{(1:k_P)}(x) - \E_{Y|X}\bar Y^P_{(1:k_P)}(x) | \leq \sqrt{\frac{(d+1)\log(n_P+n_Q) - \log(\delta/2)}{2k_P}}. \label{ada_lm1_3}  \eeq
	\beq |\bar Y^Q_{(1:k_Q)}(x) - \E_{Y|X}\bar Y^Q_{(1:k_Q)}(x) | \leq \sqrt{\frac{(d+1)\log(n_P+n_Q) - \log(\delta/2)}{2k_Q}}. \label{ada_lm1_4} \eeq
	
	Note that $\eta_{k_P,k_Q,w}(x) = w\bar Y^P_{(1:k_P)}(x) + (1-w)\bar Y^P_{(1:k_Q)}(x)$. By Cauchy-Schewatz inequality, \eqref{ada_lm1_3} and \eqref{ada_lm1_4} jointly imply
	\baligns
		|\eta_{k_P,k_Q,w}(x) - \E \eta_{k_P,k_Q,w}(x)| \leq& w|\bar Y^P_{(1:k_P)}(x) - \E_{Y|X}\bar Y^P_{(1:k_P)}(x) | \\& \quad + (1-w)|\bar Y^P_{(1:k_P)}(x) - \E_{Y|X}\bar Y^P_{(1:k_P)}(x) | \\
		\leq& \sqrt{(d+1)\log(n_P+n_Q) - \log(\delta/2)}(w(2k_P)^{-\frac{1}{2}} + (1-w)(2k_Q)^{-\frac{1}{2}}) \\
		\leq& \sqrt{(d+1)\log(n_P+n_Q) - \log(\delta/2)}\sqrt{2\p{\frac{w^2}{2k_P} + \frac{(1-w)^2}{2k_Q}}} \\
		=& \sqrt{\bigg((d+1)\log(n_P+n_Q) - \log(\delta/2)\bigg)\p{\frac{w^2}{k_P} + \frac{(1-w)^2}{k_Q}}}.
	\ealigns
	
	Therefore the above bound holds for all $k_P, k_Q, w$ with probability at least $1-\delta$.

\end{document}